\documentclass[12pt]{amsart}
\usepackage{color,enumerate,wrapfig,amssymb}
\usepackage{pxfonts}
\usepackage{graphicx}
\usepackage{subfigure}
\usepackage{eucal}

\newcommand{\R}{{\mathbb R}}

\newcommand{\N}{{\mathbb N}}

\newcommand{\osc}{\operatornamewithlimits{osc}}

\newcommand{\supp}{\operatorname{supp}}

\newtheorem{thm}{Theorem}[section]
\newtheorem{lem}[thm]{Lemma}
\newtheorem{example}[thm]{Example}
\newtheorem{prop}[thm]{Proposition}
\newtheorem{cor}[thm]{Corollary}
\newtheorem{definition}[thm]{Definition}

\newtheorem{remark}[thm]{Remark}

\title[Non-$C^2$ Hamilton-Jacobi Obstacles]{Double Obstacle Problems with obstacles given by non-$C^2$ Hamilton-Jacobi equations}

\author[John Andersson]{John Andersson}
\address{Mathematics Institute, University of Warwick Coventry CV4 7AL, UK}
\email{J.E.Andersson@warwick.ac.uk}

\author[Henrik Shahgholian]{Henrik Shahgholian}
\address{Department of Mathematics, Royal Institute of Technology,
  100~44  Stockholm, Sweden}
\email{henriksh@kth.se}

\author[Georg Weiss]{Georg S. Weiss}
\address{Department of Mathematics, Heinrich Heine University, 40225 D\"usseldorf}
\email{weiss@math.uni-duesseldorf.de}

\thanks{H. Shahgholian has been supported in part by Swedish Research Council. 
 Both J. Andersson and G.S. Weiss thank G\"oran Gustafsson's Foundations for visiting appointments to The Royal Institute of Technology}

\begin{document}

\begin{abstract}
We prove optimal regularity for the double obstacle problem when  obstacles
are given by  solutions to  Hamilton-Jacobi equations that are not $C^2$.

When the Hamilton-Jacobi equation is not $C^2$ then the standard
Bernstein technique fails and we loose the usual semi-concavity
estimates. Using a non-homogeneous scaling (different speed in different directions) 
we develop a new pointwise regularity theory for Hamilton-Jacobi
equations at points where the solution touches the obstacle.

A consequence of our result is that  $C^1$-solutions to the
Hamilton-Jacobi equation 
$$
\pm |\nabla h-a(x)|^2=\pm 1  \textrm{ in } B_1, \qquad 
h=f  \textrm{ on } \partial B_1,
$$
are in fact $C^{1,\alpha/2}$ provided that $a \in C^\alpha$. This result is optimal and  to the authors' best knowledge  new.
\end{abstract}

\maketitle

\tableofcontents

\section{Introduction.}

\subsection{Background.}

The classical torsion problem modelling the elastic-plastic torsion
of a bar can be formulated as follows:
$$
\textrm{Minimize }\int_{\Omega}\big( |\nabla u|^2 -u \big) \> dx
$$
in the convex set $K=\{u\in W^{1,2}_0(\Omega);\; |\nabla u|^2\le 1\}$ (see
\cite{T1}, \cite{T2}, \cite{T3}, \cite{CF}). Brezis and Sibony \cite{BS} showed that this
problem is equivalent to minimizing the above energy in the set
$\tilde{K}=\{u\in W^{1,2}_0(\Omega);\; u \le d(x) \}$ where $d(x)$ is
the distance function to the boundary.

More generally, one can show that minimizing the Dirichlet energy
in $L=\{u\in W^{1,2}(\Omega);\; |\nabla u|^2\le 1 \textrm{ and } u= f \textrm{ on } \partial \Omega\}$ is equivalent to minimizing the Dirichlet energy
in the set $\tilde{L}=\{u\in W^{1,2}(\Omega);\; h^-\le u \le h^+\}$
where $h^\pm$ solves the Hamilton-Jacobi equations
$$
\begin{array}{ll}
\pm|\nabla h^\pm|^2=\pm 1 & \textrm{in } \Omega, \\
h^\pm =f & \textrm{on } \partial \Omega,
\end{array}
$$
provided that $L$ and $\tilde{L}$ are non empty, which incidentally are equivalent
conditions (see \cite{BS} for a proof of this equivalence).

The regularity theory for a minimizer of the Dirichlet energy in $\tilde{L}$
is quite straightforward.  Indeed, one may approximate
$h^+$ by the solution $h^+_\epsilon$ to
$$
\begin{array}{ll}
-\epsilon \Delta h^+_\epsilon +|\nabla h^+_\epsilon|^2=1 &
\textrm{in } \Omega ,\\
h^+_\epsilon =f & \textrm{on } \partial \Omega,
\end{array}
$$
whence the Bernstein technique (see \cite{L} or Lemma \ref{OneSidedEstimates} below)
gives an $\epsilon$-independent estimate on the second derivatives
from above. In particular, letting $\epsilon\to 0$ we may deduce that
$h^+$ is semi-concave, or equivalently that the distributional second
derivatives of
$h^+$ are bounded from above. Similarly, the second order distributional
derivatives of $h^-$ are bounded from below.
From here, standard regularity theory for the obstacle
problem (such as developed in \cite{C}, modified slightly in the present paper 
to suit our purposes) implies that $u\in C^{1,1}_{\rm loc}$.

The important step in the above proof is to deduce one-sided estimates
on the distributional second derivatives of $h^\pm$. More generally, if
$u$ is a  minimizer to the Dirichlet energy in $\tilde{L}$ with
$h^+$ being a solution to a Hamilton-Jacobi equation $F(x,h^+,\nabla h^+)=1$
($h^-$ solves $-F(x,h^-,\nabla h^-)=-1$) for an $F\in C^2$ satisfying some
structural assumptions, then
the same technique yields $u\in C^{1,1}$. The important step again is the
Bernstein technique where we need to differentiate $F$ twice.
See for instance \cite{ChoeShim} or \cite{Jensen} for variational problems of 
this type or \cite{Evans1} with errata \cite{Evans2} for a fully nonlinear
gradient constrained problem.

When $F\notin C^2$, the above outline of a regularity proof fails.
The existing regularity theory for Hamilton-Jacobi
equations is not sufficiently strong for that purpose. 
There are  cases where certain one-sided estimates have been deduced, see \cite{CS} and \cite{S}. 
For example, in \cite{CS} and \cite{S} it is
assumed that $F\in C^{0,1}$ (besides standard structural assumptions)
and it is shown that $h^+(x+x^0)\le h^+(x^0)+ p\cdot x+ |x|\sigma(|x|)$
for any $p$ in the super-differential of $h^+$ at $x^0$ where $\sigma$ is
some one-sided modulus of continuity over which we have no control.

\subsection{Main Result and Ideas}

The objective of this paper is  to introduce new techniques handling regularity questions for  Hamilton-Jacobi 
equations below the $C^2$-threshold. It is noteworthy that the minimization problem in class $L$
with Hamilton-Jacobi  equations below $C^2$-threshold has applications in micro-magnetics.
In view of the equivalence of the minimization problem in class $L$ and $\tilde L$
for the special case of gradient constraint $|\nabla u|\le 1$ (see Background),
we study in the present paper the problem in $\tilde L$ with variable coefficients that
are not Lipschitz. 

Our main results are the following theorem and its corollary.

\begin{thm}\label{obstacle} [Main Theorem] 
Let $u$ be a minimizer of the Dirichlet energy
$$
\int_{B_1}|\nabla u|^2
$$
in the set $\{u\in W^{1,2};\; h^-\le u\le h^+,\; u=f\in C^{\alpha}(\partial B_1) \textrm{ on }
\partial B_1\}$, where $h^\pm$ are viscosity solutions to
\begin{equation}\label{HJ-equation}
\begin{array}{ll}
\pm |\nabla h^\pm-a(x)|^2=\pm 1 & \textrm{in } B_1, \\
h^\pm=f & \textrm{on } \partial B_1,
\end{array}
\end{equation}
with $[a]_{C^\alpha(\overline{B_1})}= A\leq \tilde C$ for some $\tilde C < {+\infty}$.
Assume furthermore that $u(0)=h^+(0)$. Then
\begin{equation}
\begin{array}{ll}
\osc_{x\in B_r(0)}\big(u(x)-u(0)-\nabla u(0)\cdot x\big)\le C(\alpha)\sqrt{A}r^{1+\frac{\alpha}{2}} &
\textrm{if } r\le A^{\frac{1}{2-\alpha}}, \\
\osc_{x\in B_r(0)}\big( u(x)-u(0)-\nabla u(0)\cdot x\big)\le C r^{2} & \textrm{if } r> A^{\frac{1}{2-\alpha}}.
\end{array}
\end{equation}
\end{thm}
A surprising consequence which may be of immediate interest to the regularity theory for
Hamilton Jacobi equations 
is the following corollary.

 \begin{cor}\label{main-cor}
Let $h$ be a $C^1$ solution to  $|\nabla h-a|^2=1$. Then  
 $h \in C^{1,\alpha/2}$, provided that $a\in C^\alpha$. 
 \end{cor}
 
 \proof
 Since $h\in C^1$, we have, by uniqueness of solutions to HJ-equations, that if $h^\pm$ solves  equations 
 $$
 \begin{array}{ll}
 \pm |\nabla h^\pm-a(x)|^2=\pm 1 & \textrm{in } B_1, \\
 h^\pm=h & \textrm{on } \partial B_1,
 \end{array} $$ 
then $h^+=h^-=h$. In particular the set 
$$
K=\{u \in W^{1,2}; \, h^- \leq u \leq h^+\}= \{h\}.
$$ 
Therefore, $h$ is in a trivial way a minimizer of $\int_{B_1} |\nabla u|^2$ in $K$. By Theorem \ref{obstacle} it follows that $h\in C^{1,\alpha/2}$.
 \qed

The function  $F(x,p)=|p-a(x)|^2$ is related  to a HJ equation that arises in  micro-magnetics, and hence   our choice  is not completely arbitrary  (see \cite{AKM}). 
 Our method is quite robust, and as such it seems plausible to adapt it to a wide class of Hamilton-Jacobi equations.

The main difficulty, as indicated above, is to develop a
strong enough regularity theory for Hamilton-Jacobi equations. In this article,
we will not treat Hamilton-Jacobi equations in their full generality. Instead we will
stay within the confines of the obstacle problem. This has
one great advantage: it is easy to see that $u\in C^{1,\beta}$ for some
$\beta>0$ (Lemma \ref{PrelimObstReg} and Proposition \ref{FirstRegObst}). 
We will therefore get a
one-sided estimate of $h^+$ from below at all points where $u=h^+$. That means that the set $\{ u=h^+\}$ 
does not intersect the singular set of $h^\pm$; by the singular set of $h^+$ 
we mean the set where $h^\pm$ are not differentiable in the classical sense.

The regularity theory for $h^+$ is deduced by {\em inhomogeneous scaling}.
There is a slight complication to apply this method to
Hamilton-Jacobi equations of our type. In particular, even though $h^+$
satisfies an elliptic Hamilton-Jacobi equation it scales parabolically
and the blow-up limit will solve a parabolic equation. Let us denote
$$
h_j(x)=\frac{h^+(r_j x',r_j^{1-\beta}x_n)}{r_j^{1+\beta}}.
$$
Then if $h_j\to h_0$ as $r_j\to 0$, the function $h_0$ heuristically
solves an equation of the form $|\tilde\nabla h_0|^2-2\partial_nh_0=0$
where $\tilde\nabla=(\partial_1,...,\partial_{n-1},0)$.

Our first goal (Proposition \ref{nondeg}) is to show that the $h_j$ defined
above is indeed bounded. It is here that we use the assumption 
$u(0)=h^+(0)$, which gives a one-sided estimate from below on $h^+$.

Once that is proved we can use the regularity theory for parabolic
Hamilton-Jacobi equations to deduce that $h^+$ satisfies
somewhat better
one-sided estimates in the $x'$ directions, for example
$h^+(x',0)\le h^+(0)+p\cdot x'+C|x'|^{1+\beta+\epsilon}$ where
$\epsilon=\epsilon(\beta,\alpha)$ (see Proposition \ref{orthogonal}).

The minimizer $u$ on the other hand scales elliptically, so
if $$u(r_jx)/\sup_{B_{r_j}}|u|$$ converges to $u_0$ as $r_j\to 0$,
then $u_0$ is a solution to an obstacle problem with
obstacle $h_0=\lim_{j\to \infty}h^+(r_j x)/\sup_{B_{r_j}}|u|$. This fact will be
used in the proof of our main theorem, Theorem \ref{obstacle}, in order to show that
$\osc_{B_r}|u(x)-\nabla u(0)\cdot x|\le C|x|^{1+\beta+\epsilon}$. In particular
we gain an $\epsilon$ in the regularity of $u$. By carefully
keeping track of all the constants we see that this can be iterated
to obtain $\gamma$-independent $C^{1,\gamma}$
estimates of $u$ for all $\gamma<\alpha/2$. 
It follows that $u\in C^{1,\alpha/2}$.

The following example pointed out to us by Stefan M\"uller 
shows that this is indeed the optimal regularity:

\begin{example}\label{sm}
If $a(x)=|x_2|^\alpha e_1$, then
$u(x_1,x_2)=x_1+b(x_2)$ is a solution to $|\nabla u-a(x)|^2=1$ for
$$
b(x_2)=\int_0^{x_2}\sqrt{2|t|^\alpha-|t|^{2\alpha}}\> dt.
$$
is a solution to $|\nabla u-a(x)|^2=1$. Here $a\in C^{\alpha}$ and
$u\in C^{1,\alpha/2}$, but $u\notin C^{1,\beta}$ for any $\beta> \alpha/2$.
\end{example}

The plan of the paper is as follows. In Section \ref{doubleObstacleSection}  
we deduce an abstract 
regularity result for solutions to the double obstacle problem, which
we will need later. In the subsequent two sections we show that $h^+$
remains bounded under parabolic scaling and that $h^+$ satisfies better
one-sided estimates in the $x'$ directions. In the final section
we prove our main result that the minimizer $u\in C^{1,\alpha/2}$.
Finally we have included a long appendix to remind the reader
of some of the theory of viscosity solutions for Hamilton-Jacobi equations.
In the Appendix we also deduce simple $C^{1,\beta}$-estimates for
the solution, which will serve as our starting regularity in the 
strategy described above.

\textbf{Acknowledgment:} We would like to thank Stefan M\"uller for providing us with 
Example \ref{sm}.

\subsection{Notation} 
We denote the Euclidean ball $B_r(x) = \{ y\in \R^n; |y-x|<r\}$, and we denote
$\omega_n := \mathcal{L}^n(B_1)$; in the case
that the center is not specified it is assumed to be the origin.
When $v\in \mathbb{R}^n$ is a vector we will denote the first $n-1$ coordinates by $v':= (v_1,v_2,...,v_{n-1})$. Similarly, we will use
$\tilde\nabla:=(\partial_1,\partial_2,...\partial_{n-1},0)$. 
Here $\partial_i=\frac{\partial }{\partial x_i}$ denotes differentiation
with respect to the $x_i$ variable. We will often denote differentiation
by a subscript $\partial_i u =: u_i$. The unit vector in the $i-$th
coordinate direction will be denoted by $e_i$.
We will also use the seminorm
$$
[a]_\alpha=[a]_{C^{\alpha}(\bar \Omega)}:= \sup_{x,y\in \bar \Omega}\frac{|a(x)-a(y)|}{|x-y|^\alpha}.
$$
Finally, different instances of the letter $C$ may mean different constants
even within one proof or one set of inequalities.

\section{Regularity for Double Obstacle Problems}\label{doubleObstacleSection}

In this section we will prove regularity for the two-sided obstacle 
problem in the context needed later.

It will be convenient to define a class of solutions to double obstacle problems
and to fix some notation. We therefore state the following two definitions before we
state and prove the main result of this section.
\begin{definition}
For any continuous function $u$ we define the super-differ\-en\-tial 
(sub-differential) of $u$ at the point $x^0\in \textrm{Domain}(u)$
as follows:
$$
S^+(u,x^0)=\big\{ p\in \mathbb{R}^n;\; \sup_{B_r}\big(u(x+x^0)-u(x^0)-p\cdot (x-x^0)\big)\le o(r)  \big\}
$$
$$
\Big( S^-(u,x^0)=\big\{ p\in \mathbb{R}^n;\; \inf_{B_r}\big(u(x+x^0)-u(x^0)-p\cdot (x-x^0)\big)\ge -o(r)  \big\}\Big).
$$
\end{definition}

\begin{definition}\label{defclassC}
We define $\mathcal{C}(R,\sigma,h^\pm)$ as the set of 
local minimizers $u$ to the Dirichlet energy
$$
\int_{B_R}|\nabla u|^2
$$
in the set $K=\{ u\in W^{1,2}_{\rm loc}(B_R); h^-\le u\le h^+\}$;
here $0< R\in \mathbb{R}$, $\sigma(r)$
is a ``one-sided modulus of continuity'' of $u$ and $h^\pm\in C^0(B_R)$ such that $h^+\ge h^-$.

Furthermore we require that $u$ satisfies the following estimates for each point 
$x^0\in \Lambda(u,R):= \big(\{u=h^+\}\cup \{u=h^-\} \big)\cap B_R$:
\begin{itemize}
    \item[(i)] when $u(x^0)=h^+(x^0)$, then we have the bound from above
$$
\sup_{B_r(x^0)}\inf_{p_{x^0}\in \mathbb{R}^n}\big( u(x)-p_{x^0}\cdot(x-x^0)- u(x^0)\big) \le \sigma(r)
\textrm{ for all } r>0.
$$
\item[(ii)] when $u(x^0)=h^-(x^0)$, then we have the bound from below
$$
\inf_{B_r(x^0)}\sup_{p_{x^0}\in \mathbb{R}^n}\big( u(x)-p_{x^0}\cdot(x-x^0)-u(x^0)\big)\ge -\sigma(r)
\textrm{ for all }r>0.
$$
\end{itemize}
\end{definition}

\begin{remark}
It is important that by (i) and (ii), $u$ inherits its ``one-sided modulus of continuity'' from
that of $h^\pm$.
\end{remark}

\begin{prop}\label{justneedwhatevername}
Let $u\in \mathcal{C}(R,\sigma,h^\pm)$
where $\sigma(r)$ is a one-sided modulus of continuity satisfying the doubling
condition
\begin{equation}\label{somedoubling}
\frac{\sigma(Sr)}{\sigma(r)}\le F(S)  \textrm{ for each }S\ge 1 \textrm{ and all }r>0.
\end{equation}
for each $S\ge 1$.
Then if $x^0\in\Lambda(u,R/2)$ it follows that
$$
\osc_{B_r(x^0)}\inf_{p_{x^0}\in \mathbb{R}^n}\big(u-p_{x^0}\cdot (x-x^0)\big)
\le C(n,F) \sigma(r).
$$
\end{prop}
\textbf{Remark:} It is of no importance that the set $K$ has the particular
form $K=\{u\in W^{1,2};\; h^-\le u\le h^+\}$. We also do not need that $h^\pm$
solves any particular Hamilton-Jacobi equation in this Proposition. The important
thing is that we have one-sided estimates at points where $u$ is not harmonic.
Similarly the boundary values $f$ are of minor importance.

If we consider the function $F$ in (\ref{somedoubling}) to be  $F(s)=s^\alpha$
then the Proposition slightly generalizes known regularity results for double obstacle problems.

\textsl{Proof of Proposition \ref{justneedwhatevername}:}
We argue by contradiction and assume that there is a sequence
$u_j$ of minimizers as in the proposition, with $\sigma=\sigma_j$ ($\sigma_j$ 
corresponds to a fixed $F$) and points $x^j\in \Lambda(u,R)$ 
(without loss of generality we are going to assume that $x^j=0$), $p_{j}\in S^+(u,x)$ or
$p_j\in S^-(u,0)$
as well as $r_j>0$ such that
\begin{equation}\label{ijn}
\osc_{B_{1}}\inf_{p_j\in \mathbb{R}^n}
\frac{u_j(r_j x)-r_j p_j\cdot x}{j\sigma_j(r_j)}=1.
\end{equation}
We may also assume that for large $j$, $r_j$ is the largest such $r$ corresponding to $u^j$, that is for $p\in S^+(u,0)$ or $p\in S^-(u,0)$
\begin{equation}\label{answer7and8}
\osc_{B_{1}} \frac{u_j(r x)-r p\cdot x}{j\sigma_j(r)}\le 1
\end{equation}
for $r\ge r_j$.

Next, we define
$$
v^j(x)=\frac{u_j(r_j x)-r_j p_j\cdot x-u_j(0)}{j\sigma(r_j)},
$$
where $p_j$ is the vector in equation (\ref{ijn}).
Then 
$v^j\in \mathcal{C}(R,\sigma(r_j\cdot )/(j\sigma_j(r_j)), \tilde{h}^\pm_j)$ where
$$
\tilde{h}^\pm_j(x)=\frac{h^{\pm}_j(r_j x)-r_j p_j\cdot x-u_j(0)}{j\sigma(r_j)}.
$$
In particular
$\Delta v^j=0$ in some set $\Omega_j= B_{R/r_j}\setminus \Lambda(v^j,R/r_j)$, and for every point
$x^0 \in \Lambda(v^j,R/(2r_j))$ we have the following: if $v^j(x^0)=\tilde{h}_j^+(x^0)$ then
(cf. condition (i) in Definition \ref{defclassC})
\begin{equation}\label{concavityestimate}
\sup_{B_r(x^0)}\inf_{p\in \mathbb{R}^n}\big(v^j(x)-p\cdot (x-x^0)-v^j(x^0) \big)\le
\frac{\sigma_j(rr_j)}{j\sigma(r_j)}\le \frac{F(r)}{j}\to 0.
\end{equation}
Moreover, if $v^j(x^0)=\tilde{h}_j^-(x^0)$ then
it follows from condition  (ii)  in Definition \ref{defclassC} that
\begin{equation}\label{convexityestimate}
\sup_{p}\inf_{B_r(x^0)}\sup_{p\in \mathbb{R}^n}\big(v^j(x)-p\cdot (x-x^0)-v^j(x^0) \big)\ge
-\frac{\sigma_j(rr_j)}{j\sigma(r_j)}\le \frac{F(r)}{j}\to 0.
\end{equation}
Let us write $\Lambda(v^j,R/r_j)=\Lambda^+(v^j,R/r_j)\cup \Lambda^-(v^j,R/r_j)$ 
with (\ref{concavityestimate})
satisfied for points in $\Lambda^+$ and (\ref{convexityestimate}) satisfied
in $\Lambda^-$.

For some subsequence $v^j\to v^0$  in $L^p_{\textrm{loc}}$, where 
$$
v^0\in \mathcal{C}(t, 0, \tilde{h}^\pm_0)
$$ 
for all $t<{+\infty}$;
here we have used (\ref{concavityestimate}) and (\ref{convexityestimate}) which imply
that
$$
\sigma(r_j r)/(j\sigma_j(r_j))\to 0 \quad \textrm{ as }j\to \infty,
$$
as well as the notation $\tilde{h}^{\pm}_0=\lim_{j\to \infty}\tilde{h}^\pm$.
We will actually show later that this convergence is locally uniform.
Furthermore, by equation (\ref{concavityestimate})
we have that $S^+(v^0, x^0)\ne \emptyset$ for all $x^0\in \Lambda^+(v^0, {+\infty})$ and
$$
v^0(x)\le \inf_{x^0\in \Lambda^+}\inf_{p^+_{x^0}\in S^+(v^0,x^0)}\big(  v^0(x^0)+p^+_{x^0}\cdot(x-x^0)\big),
$$
and by equation (\ref{convexityestimate}) we have
that $S^-(v^0, x^0)\ne \emptyset$ for all $x^0\in \Lambda^-(v^0, {+\infty})$ and that
$$
v^0(x)\ge \sup_{x^0\in \Lambda^-}\sup_{p^-_{x^0}\in S^-(v^0,x^0)}\big( v^0(x^0)+p^-_{x^0}\cdot (x-x^0)  \big).
$$

That is, $v^0$ solves the double obstacle problem with a concave upper obstacle
and a convex lower obstacle. It follows that $\Delta v^0=0$ in $\mathbb{R}^n$.
By our assumption that the origin is in $\Lambda(v^0,{+\infty})$, say 
$0\in \Lambda^+(v^0,{+\infty})$, it follows that for $p\in S^+(v^0,0)$
$$
\sup_{B_r}\big( v^0(x)-p\cdot x \big) = 0
$$
for all $r<{+\infty}$.

It follows from Liouville's Theorem that $v^0$ is a linear function, but $v^0(0)=0$ 
and by construction (cf. equation (\ref{ijn})) 
\begin{equation}\label{anhour}
\osc_{B_1}\inf_{p\in \mathbb{R}^n}\big( v^j(x)-p\cdot x\big)=1.
\end{equation}
Equation (\ref{anhour}) is a contradiction to $v^0$ being linear, provided that we can show that 
$v^j\to v^0$ uniformly.

In order to show that $v^j\to v^0$ uniformly we notice that the obstacle functions
$$
f_j^+:= \inf_{x^0\in \Lambda^+_j}\inf_{p_{x^0}\in S^+(v^j,x^0)} \Big(v^j(x^0)+p_{x^0}\cdot (x-x^0)\Big)+
\frac{\sigma(|x|r_j)}{j\sigma(r_j)}
$$
converge either locally uniformly to $+\infty$ or that locally uniformly
$\lim_{j\to \infty}f_j^+=f^+$ where $f^+$ is an affine linear function. Similarly,
$$
f_j^-:= \sup_{x^0\in \Lambda^-_j}\sup_{p_{x^0}\in S^-(v^j,x^0)} \Big( v^j(x^0)+p_{x^0}\cdot (x-x^0)
\Big)-\frac{\sigma(|x|r_j)}{j\sigma(r_j)}
$$
will converge either locally uniformly to $-\infty$ or  locally 
uniformly
$\lim_{j\to \infty}f_j^-=f^-$ where $f^-$ is an affine linear function.

Since $f^-_j\le u^j\le f^+_j$ we have $f^-_j\le f^+_j$.

We may distinguish three cases:

1) $v^j(\hat{x}^j)=f_j^+(\hat{x}^j)$ for some bounded sequence of points $\hat{x}^j$ and
$v^j(\hat{y}^j)=f^-_j(\hat{y}^j)$ for some bounded sequence of points $\hat{y}^j$,

2) $v^j(\hat{x}^j)=f_j^+(\hat{x}^j)$ for some bounded sequence of points $\hat{x}^j$ and
$v^j>f^-_j$ in $B_{R_j}$ with $R_j\to {+\infty}$, or $v^j(\hat{y}^j)=f_j^-(\hat{y}^j)$ for some bounded sequence of points $\hat{y}^j$ and
$v^j<f^+_j$ in $B_{R_j}$ with $R_j\to {+\infty}$,

3) $f^-_j< v^j< f^+_j$ in $B_{R_j}$ for a sequence $R_j\to {+\infty}$.

In the first case it is easy to see that $|f^+_j-f^-_j|\to 0$ locally uniformly:
if not then $f^-_j\to a+p\cdot x$ and $f^+_j\to b+p\cdot x$ with $b>a$.
Also $v^j\to c+p\cdot x$ in $L^2_{\rm loc}$ with $b\ge c\ge a$. We may assume
that $p=0$. Now only one of the inequalities $c>a+(b-a)/3$ or $b-(b-a)/3>c$ may hold; let us for definiteness assume
that $b-(b-a)/3>c$. Then since $v^j\to c$ in $L^{2}_{\rm loc}$ there is for 
$T:=2\sup_{j}\max(|\hat{x}^j|,|\hat{y}^j|)$ and each
$\epsilon>0$ a $j_{\epsilon}<+\infty$ such that
\begin{equation}\label{u+issmall}
\frac{1}{(2T)^n}\int_{B_{2T}}|(v^j-c)^+|^2< \epsilon
\end{equation}
for $j>j_\epsilon$.
Also, as $f_j^\pm$
converges uniformly, $f^-_j< c+\epsilon$ for all $j>j_\epsilon$
(provided that $j_\epsilon$ has been chosen large enough). In particular
$\{v^j> c+\epsilon\}\cap \Lambda_j^-=\emptyset$ when $j>j_\epsilon$.
That means that $v^j$ may touch only the upper obstacle in
$\{v^j>c+\epsilon\}$. Thus $v^j$ is subharmonic in $\{v^j>c+\epsilon\}$.
From (\ref{u+issmall}) we deduce that there exists $\tilde{T}\in (3T/2,2T)$
such that
\begin{equation}\label{lifeisshit}
\frac{1}{\tilde{T}^{n-1}}\int_{\partial B_{\tilde{T}}}(v^j-c)^+\le C\epsilon
\end{equation}
for some fixed $C\ge 1$. Now we define $w^j$ by the Poisson integral as 
\begin{equation}\label{poissonwj}
w^j(x)=
\frac{\tilde{T}^2-|x|^2}{n\omega_n \tilde{T}}\int_{\partial B_{\tilde{T}}}\frac{\max(v^j,c+C\epsilon)}{|x-y|^n}dy.
\end{equation}
Then $\Delta w^j=0$ in $B_{\tilde{T}}$ and $w^j\ge v^j$ on $\partial B_{\tilde{T}}$. 
Since $\Delta v^j\ge 0$ in $\{v^j>c+\epsilon\}$ and $w^j\ge c+C\epsilon$ it
follows that $\max(v^j, w^j)$ is subharmonic. But $w^j\ge v^j$ on $\partial B_{\tilde{T}}$ and $\Delta w^j=0$
so $\max(v^j,w^j)\le w^j\le \max(v^j,w^j)$, where the first inequality follows by comparison and the second
is trivial. Thus $w^j\ge v^j$. Moreover we know from (\ref{poissonwj}) and (\ref{lifeisshit}) that
$w^j\le c+C\epsilon$ in $B_{\tilde{T}}$. If $\epsilon$ is small enough we may deduce that in $B_{\tilde{T}}$,
$$
v^j\le c+C\epsilon <b-\frac{b-a}{3}+C\epsilon< b-\frac{b-a}{6}<f^+_j.
$$
But that contradicts the condition in 1) that $v^j$ touches the upper obstacle.

In case 2) we may argue similarly as in case 1) to show that for each $T$,
$f^+_j\ge v^j> f^+_j-\epsilon$ in $B_T$ if $j$ is
large enough, or
$f^-_j\le v^j< f^-_j+\epsilon$ in $B_T$ if $j$ is
large enough, respectively.

In case 3) we have that for each $R<+\infty$, $v^j$ is harmonic in $B_{R}$ 
for large $j$, which implies
uniform convergence in $C^2_{\rm loc}$. The Proposition follows. \qed

Using Lemma \ref{PrelimObstReg} we can derive a simple Corollary that states some
preliminary regularity for solutions to the obstacle problem.

\begin{cor}\label{FirstRegObst}
Let $u$ be a minimizer to the two-sided obstacle problem as above,
i.e. $u$ is the minimizer of
$$
\int_{B_R} |\nabla u|^2
$$
in the set $\{u\in W^{1,2}; \,   h^-\le u\le h^+\}$ and
$h^\pm$ are solutions to $\pm|\nabla h^\pm-a|^2= \pm 1$ where
$[a]_{\alpha}\le A$.
Then if $x^0\in B_{R/2}$ and $u(x^0)=h^+(x^0)$ (or $u(x^0)=h^-(x^0)$) and $p_{x^0}$ is in the 
super(sub)-differential of $u$ at $x^0$, it follows that
\begin{equation}\label{cromangon}
\osc_{B_r} \left(u(x+x^0)-p_{x^0}\cdot (x-x^0)\right)\le C(n) A^{1/(2+\alpha)}r^{1+\alpha/(2+\alpha)}
\end{equation}
for $r\le CA^{1/2}R^{\frac{2+\alpha}{2}}$ and that
\begin{equation}\label{cromangon2}
\osc_{B_r} \left(u(x+x^0)-p_{x^0}\cdot (x-x^0)\right)\le \frac{C}{R}r^2
\end{equation}
for $r\ge CA^{1/2}R^{\frac{2+\alpha}{2}}$.

In particular, $u\in C^{1,\frac{\alpha}{2+\alpha}}(\overline{B_{R/2}})$.
\end{cor}
\textsl{Proof:} Using Lemma \ref{PrelimObstReg} we see that
$u$ satisfies the assumptions of Proposition \ref{justneedwhatevername}
with $\sigma$ defined by
$$
\sigma(r)= C A^{1/(2+\alpha)}r^{1+\alpha/(2+\alpha)}
$$
for $r\le CA^{1/2}R^{(2+\alpha)/2}$,
and
$$
\sigma(r)= C \frac{r^2}{R},
$$
for $r\ge CA^{1/2}R^{(2+\alpha)/2}$. 

It follows then by a standard method that
$u\in C^{1,\frac{\alpha}{2+\alpha}}$: from (\ref{cromangon})
it follows that $\nabla u(x^0)$ is well defined. Whitney's extension theorem 
implies that we can extend $u\big|_{\supp(\Delta u)}$ to a function $\hat{u}$
defined in $B_{R/2}$ where $\nabla \hat{u}$ is $C^{\alpha}$ with modulus of continuity
$$
\sigma(r)=\left\{ \begin{array}{ll}
C(n)A^{1/(2+\alpha)}r^{\alpha/(2+\alpha)} & \textrm{for }r\le CA^{1/2}R^{\frac{2+\alpha}{2}}, \\
\frac{C}{R}r & \textrm{for } r>CA^{1/2}R^{\frac{2+\alpha}{2}}.
\end{array}\right.
$$
In particular if $x,y\in \supp(\Delta u)$ then 
$|\nabla u(x)-\nabla u(y)|\le \sigma(|x-y|)$.

Also notice that if $\textrm{dist}(x^0, \supp(\Delta u))=|x^0-y|=\kappa$ 
where $y\in \supp(\Delta u)$ then 
$\Delta u=0$ in $B_{\kappa}(x^0)$ and the supremum of $|u(x)-u(y)-\nabla u(y)\cdot (x-y)|$ 
in $B_{\kappa}(x^0)$ can be estimated by $\sigma(2\kappa)$. Standard regularity theory for harmonic
functions implies that $|\nabla u(x^0)-\nabla u(y)|\le C\sigma(2\kappa)$. In particular
it follows that if $x\in B_{R/2}\setminus \supp(\Delta u)$ and 
$y\in B_{R/2}\cap \supp(\Delta u)$ then $|\nabla u(x)-\nabla u(y)|\le C\sigma(2|x-y|)$.

If both $y,z\in B_{R/2}\setminus \supp(\Delta u)$ and 
$$
|y-z|<\frac{1}{2}
\max(\textrm{dist}(y,\supp(\Delta u)),\textrm{dist}(z,\supp(\Delta u)))=\kappa,
$$
then $u$ is harmonic in $B_\kappa(y)$ or in $B_{\kappa}(z)$. For the sake of definiteness
we will assume that $u$
is harmonic in $B_\kappa(y)$. 
Also
$\textrm{sup}_{B_{\kappa}(y)}|u(x)-u(y)-\nabla u(y)\cdot (x-y)|\le \sigma(2\kappa)$. Standard 
regularity theory for harmonic functions implies that 
$$
|\nabla u(y)-\nabla u(z)|\le C\sigma(2\kappa)|y-z|/\kappa\le C\sigma(|y-z|).
$$
Finally, if $y,z\in B_{R/2}\setminus \supp(\Delta u)$ and
$$
|y-z|\ge \frac{1}{2}
\max(\textrm{dist}(y,\supp(\Delta u)),\textrm{dist}(z,\supp(\Delta u)))=\kappa,
$$
then we may combine the first two estimates and deduce that 
$$
|\nabla u(y)-\nabla u(z)|\le C\sigma(|y-z|).
$$
The Corollary follows.
\qed

%%%%%%%%%%%%%%%%%%%%%%%%%%%%%%%%%%%%%%

%%%%%%%%%%%%%%%%%%%%%%%%%%%%%%%%%%%%%%%%%%

\section{Improved Regularity in the Non-Degenerate Direction}

In this section we will discuss the behavior of a viscosity solution
in the non-degenerate direction. To motivate this discussion
and our terminology, let us consider a simple elliptic rescaling
of $|\nabla h^+-a|^2=1$: if $h^+$ is differentiable at the origin,
$h^+(0)=|\nabla h^+(0)|=0$, $a(x)=e_n+b(x), |b(0)|=0$ and
$$
h_j(x)=\frac{h^+(r_jx)}{\sup_{B_{r_j}}|h^+|}
$$
converges to $h_0$ as $r_j\to 0$. If we rescale we see that 
\begin{align*}
0=&|\nabla h_j|^2\frac{\sup_{B_{r_j}}|h^+|}{r_j}-2\frac{\partial h_j}{\partial x_n}-
2b(r_j x)\cdot \nabla h_j\\
&\; +|b(r_jx)|^2\frac{r_j}{\sup_{B_{r_j}}|h^+|}+2b_n(r_jx)
\frac{r_j}{\sup_{B_{r_j}}|h^+|},
\end{align*}
and ---at least heuristically---
$$
\frac{\partial h_0}{\partial x_n}=0,
$$
provided that
$$|b(r_jx)|^2\frac{r_j}{\sup_{B_{r_j}}|h^+|}+2b_n(r_jx)
\frac{r_j}{\sup_{B_{r_j}}|h^+|}\to 0, j\to\infty.$$
Therefore the Hamilton-Jacobi
equation is degenerate in all directions except one. Thus we
might expect that the oscillation of $h^+(x',x_n)$  along lines 
$\{(x',s);\; |s|\le r\}$
with fixed $x'$ should be of lower order compared to $\osc_{B_r}h^+$. In the next
Proposition we will show that is indeed the case, if we have a
good estimate on $h^+$ from below.

The proof, and even the statement, is quite long and not very easy to read.
Therefore we will try to explain the general idea before we start.

The idea of the proof is that $h^+$ scales like $h^+(rx',r^{1-\beta}x_n)/r^{1+\beta}$,
whenever $h^+(0)=|\nabla h^+(0)|=0$ (notice that if $u(0)=h^+(0)$ then $\nabla h^+$
is defined at the origin).
It is therefore a natural assumption that $h^+(rx',r^{1-\beta}x_n)/r^{1+\beta}$
will be bounded if $u\in C^{1,\beta}$ and $h^+$ has one-sided
$C^{1+\beta}$-estimates. We will prove this boundedness (or a slightly refined version
of it) by blow-up and a contradiction argument. Most of the proof
consists of technical estimates of different terms in the equation
for the rescaled  $h^+$. The idea is simple though. If the supremum of the
rescaled $h^+$ goes to infinity then  we can find another rescaling
(called $h_j$ in the Proposition)
that respects the natural scaling of the Hamilton-Jacobi equation
(Claim 1, below) and is bounded. This new rescaling of $h^+$
has the nice property that it is worse than a $C^{1,\beta}$ scaling in the
$x'$ directions. This will imply that $h_j\to h_0$ where $h_0$ is independent
of the $x'$ directions, at least at $x_n=0$ (Claim 2 in the Proposition).
The remainder of the proof consists in verifying that $h_0$ also satisfies a
good Hamilton-Jacobi equation that will allow us to conclude that
$h^0\equiv 0$ which in turn will imply a contradiction to the fact that $h_0(0)=0$
and $\sup_{B_1}|h_0|=1$.

\begin{remark}
In the sequel we may assume that $A > 0$.
The case $A=0$ may then be handled by approximating $0< A_j \to 0$.
\end{remark}

We will denote 
\begin{equation}\label{KrDef}
K(r)=\Big\{x;
\big( \frac{x'}{r}, \frac{A^{\frac{1-\beta}{2}}x_n}{r^{1-\beta}} \big)\in B_1
\Big\}.
\end{equation}

Keeping this definition in mind, let us now state the main result in this section.

\begin{prop}\label{nondeg} Let $h^+$ be a viscosity solution to
$|\nabla h^+-a|^2=1$ in $B_R$ such that $h^+(0)=0$. Assume furthermore
\begin{itemize}
\item[(a1)]
that $h^+$ satisfies one-sided $C^{1,\beta}$ estimates:
\begin{equation}\label{one}
\sup_{B_r(x^0)}\big( h^+(x) -h^+(x^0)-p\cdot (x-x^0)\big)\le CA^{\frac{1-\beta}{2}}r^{1+\beta}
\end{equation}
for every $x^0$ and every $p$ in the super-differential of $h^+$ at $x^0$, 
\item[(a2)] that
\begin{equation}
h^+(x',x_n)\le h^+(0,x_n)+p\cdot x' + CA^{\frac{1-\beta}{2}} |x'|^{1+\beta}
\end{equation}
for each fixed $x_n$,
\item[(a3)] that 
\begin{equation}\label{onesidedfrombelowinprop2}
h^+(x)\ge -C A^{\frac{1-\beta}{2}}|x|^{1+\beta},
\end{equation}
\item[(a4)]
that $[a]_{C^\alpha}\le A$ where $A$ satisfies $A \le \tilde{C}$ for some $\tilde{C}$,
\item[(a5)]
that $a(x)=e_n+b(x)$ where $|b(0)|=0$, 
\item[(a6)] and that 
$\alpha/(2+\alpha)\le \beta\le \alpha/2$.
\end{itemize}
Then, for $r\le \min(A^{1/2}R^{\frac{1}{1-\beta}},R)$,
$$
\sup_{K(r)}|h^+(x)+g(x_n)|\le CA^{\frac{1-\beta}{2}}r^{1+\beta},
$$
where 
\begin{equation}\label{g}
g(0)=g'(0)=0,\quad g'(t)=-e_n\cdot b(0',t),
\end{equation}
and $K(r)$ is as defined in (\ref{KrDef}).
\end{prop}

\begin{remark}
The condition $r\le \min(A^{1/2}R^{\frac{1}{1-\beta}},R)$
assures that $K(r)\subset B_R$ but is otherwise not used in the proof.
\end{remark}

\begin{remark}\label{ab}
The assumption (a5) may always be satisfied scaling and rotating, so that
it is not restrictive at all.
\end{remark}

\proof We will argue by contradiction and assume that there
exist sequences $h_j^+$, $r_j$, $A_j$, $b^j$ and $a^j$ satisfying the hypothesis, such that
$$
\sup_{K(r_j)}|h_j^++g_j|=jA_j^{\frac{1-\beta}{2}}r_j^{1+\beta},
$$
where $g_j$ satisfies (\ref{g}).
The proof is rather long so we will divide it into several claims. The
first one will slightly modify $K(r_j)$ so that it respects the natural scaling
of the Hamilton-Jacobi equation.

\vspace{3mm}
\noindent
\textbf{Claim 1: }{\sl For
$$
S_j(\sigma)=\sup_{B_1} \big(h_j^+(\sqrt{\sigma} r_j x', r_j^{1-\beta}A_j^{-\frac{1-\beta}{2}}x_n)+
g_j(r_j^{1-\beta}A_j^{-\frac{1-\beta}{2}}x_n)\big),
$$
there exists $\sigma_j$ such that $S_j(\sigma_j)=\sigma_jA_j^{\frac{1-\beta}{2}}r_j^{1+\beta}$
and $\sigma_j\to {+\infty}$.}

\noindent
\textsl{Proof of Claim 1:} By assumption
$S_j(1)=jA_j^{\frac{1-\beta}{2}}r_j^{1+\beta}> > A_j^{\frac{1-\beta}{2}}r_j^{1+\beta}$.
Setting $\tau_j:=1/(A_j^{1-\beta}r_j^{2\beta})$ and using the one-sided estimates for $h^+$  (inequality (\ref{one})), we also obtain that
$$
S_j(\tau_j)\le \sup_{B_{r_j^{1-\beta}A_j^{-(1-\beta)/2}}}|h_j^+|  +
\sup_{B_{r_j^{1-\beta}A_j^{-(1-\beta)/2}}}|g_j|
$$
$$\le
CA_j^{\frac{\beta^2-\beta}{2}}r_j^{1-\beta^2}+
r_j^{(1-\beta)(1+\alpha)}A_j^{1-\frac{(1-\beta)(1+\alpha)}{2}}
$$
$$\le
C\Big(\tau_j r_j^{1+\beta}A_j^{\frac{1-\beta}{2}} \Big)
A_j^{\frac{1-2\beta+\beta^2}{2}}r_j^{\beta-\beta^2}< \tau_j A_j^{(1-\beta)/2}r_j^{1+\beta},
$$
provided that $r_j$ is small enough.

By the continuity of $S_j$ there is a $\sigma_j\in (1,\tau_j)$ such that
$S_j(\sigma_j)= \sigma_j A_j^{\frac{1-\beta}{2}}r_j^{1+\beta}$.

Therefore we only need to show that $\sigma_j\to {+\infty}$. By definition we have
$$
\sigma_j=
\sup_{B_1}\frac{h_j^+(\sqrt{\sigma_j}r_jx',r_j^{1-\beta}A_j^{-(1-\beta)/2}x_n)+g_j(r_j^{1-\beta}A_j^{-(1-\beta)/2}x_n)}{A_j^{(1-\beta)/2}r_j^{1+\beta}}
$$
$$
\ge \sup_{B_1}
\frac{h_j^+(r_jx',r_j^{1-\beta}A_j^{-(1-\beta)/2}x_n)+g_j(r_j^{1-\beta}A_j^{-(1-\beta)/2}x_n)}{A_j^{(1-\beta)/2}r_j^{1+\beta}}=j,
$$
which proves Claim 1.

Let us also remark that as
$$
S_j(\sigma_j)=\sigma_j A_j^{\frac{1-\beta}{2}}r_j^{1+\beta},
$$
 we obtain from our one-sided estimates that
$$
S_j(\sigma_j)\le
C\max\big( \sup_{B_{\sqrt{\sigma_j}r_j}}h_j^+  + \sup_{B_{\sqrt{\sigma_j}r_j}} |g_j|, \sup_{B_{A_j^{-(1-\beta)/2}r_j^{1-\beta}}}h_j^+ + \sup_{B_{A_j^{-(1-\beta)/2}r_j^{1-\beta}}} |g_j|\big).
$$
Notice that, since $g'_j(t)=b(0',t)$,
$$
\sup_{B_{A_j^{-(1-\beta)/2}r_j^{1-\beta}}} |g_j| \leq
C\left( A_j^{\frac{1-\alpha + 2\beta - (1+ \alpha)\beta}{2}} r_j^{(\alpha - \beta)(1-\beta)}\right)A_j^{\frac{\beta^2 -\beta}{2}}r_j^{1-\beta^2},
$$
but the term in the parenthesis tends to zero, implying that
$$
\sup_{B_{A_j^{-(1-\beta)/2}r_j^{1-\beta}}} |g_j| \leq A_j^{\frac{\beta^2 -\beta}{2}}r_j^{1-\beta^2}.
$$
Similarly, 
$$
\sup_{B_{\sqrt{\sigma_j}r_j}} |g_j| \leq A_j\sigma_j^{(1+ \alpha)/2}r_j^{1+\alpha} \leq 
\left(\sigma_j^{\frac{\alpha - 1}{2}}A_j^{\frac{1+ \beta}{2}} r_j^{\alpha -\beta})   \right)
\sigma_jA_j^{\frac{1- \beta}{2}}r_j^{1+\beta} \leq S_j(\sigma_j).
$$
Using these estimates on $|g_j|$ we arrive at
$$
S_j(\sigma_j)\le C\max\big( \sigma_j^{\frac{1+\beta}{2}}A_j^{\frac{1-\beta}{2}}r_j^{1+\beta}, A_j^{\frac{\beta^2-\beta}{2}}r_j^{1-\beta^2}\big).
$$
Since $\sigma_j\to {+\infty}$ it follows that
$$
\sigma_j A_j^{\frac{1-\beta}{2}}r_j^{1+\beta}\le 
C A_j^{\frac{\beta^2-\beta}{2}}r_j^{1-\beta^2}.
$$
We may ---in order to simplify notation--- assume that $\sigma_j=j$ and 
rewrite the above relation as
\begin{equation}\label{relationrjandj}
r_j\le C A_j^{-\frac{1-\beta}{2\beta}}j^{-\frac{1}{\beta+\beta^2}}.
\end{equation}

As mentioned before, Claim 1 gives an estimate that respects the scaling of
$h_j^+$: if we define $h_j$ by
$$
h_j(x):=
\frac{h_j^+(\sqrt{j}r_jx',r_j^{1-\beta}A_j^{-\frac{1-\beta}{2}}x_n)+g_j(r_j^{1-\beta}A_j^{-\frac{1-\beta}{2}}x_n)}{jA_j^{\frac{1-\beta}{2}}r_j^{1+\beta}}
$$
then $\sup_{B_1}|h_j|=1$. Also
$$
\nabla h_j(x)= \frac{\tilde\nabla h_j^+(y)}{\sqrt{j}A_j^{\frac{1-\beta}{2}}r_j^\beta}+\frac{\partial_n h_j^+(y)+g'_j(y_n)}{jA_j^{1-\beta}r_j^{2\beta}}e_n
$$
where $y=\big( \sqrt{j}r_j x',r_j^{1-\beta}A_j^{-(1-\beta)/2}x_n\big)$.
Alternatively we may write
$$
\nabla h_j^+= c_j\tilde\nabla h_j +\big( c_j^2\partial_nh_j-g_j'  \big)e_n,
$$
where
$$
c_j =\sqrt{j}A_j^{\frac{1-\beta}{2}}r_j^{\beta}.
$$
Rewriting the Hamilton-Jacobi equation (\ref{HJ-equation})
in terms of our new function $h_j$ and using assumption (a5), we may deduce that
\begin{equation}\label{eq:HJforhj}
    \begin{array}{ll}
&|\tilde\nabla h_j|^2-2\partial_n h_j=\\
&\frac{2b\cdot\tilde\nabla h_j}{c_j}
+2\partial_n h_j (b\cdot e_n + g_j')
- \frac{|b + e_n g'_j|^2}{c_j^2}
-2\frac{e_n\cdot b + g_j'}{c_j^2}-c_j^2(\partial_n h_j)^2=\\
&T^j_1\cdot \tilde\nabla h_j+T^j_2+T^j_3+T^j_4-T^j_5.
\end{array}
\end{equation}

In order to use this equation we need to control the right-hand side.
First we need to control the one-sided oscillation in the $x'$ directions. 
\vspace{3mm}
\noindent
\textbf{Claim 2: }
{\sl For every $x_n$ and every $p_j$ in the super-differential of $h_j$
at $(x',x_n)$ we have for each $\tilde{R}$ and any $q_j$ in the super-differential of $h_j^+$ at
$(0,r_j^{1-\beta}A_j^{-\frac{1-\beta}{2}})$ that
$$
\sup_{B_{\tilde{R}}'}\big(h_j-h_j(0,x_n)-p_j\cdot x'\big)  \to 0 \textrm{ as } j\to\infty.
$$
In particular, if $x_n=0$ then $\osc_{B_{\tilde{R}}'}(h_j(x',0)-h_j(0)-p_j\cdot x')\to 0\textrm{ as } j\to\infty$.}

\noindent
\textsl{Proof of Claim 2:} The first statement is deduced from (a2) in the following way:
$$
\sup_{B_{\tilde{R}}'}\big(h_j-h_j(0,x_n)-p_j\cdot x'\big)\le \sup_{B'_{\sqrt{j}r_j\tilde{R}}} 
\frac{h_j^+-h_j^+(0,x_n)+q_j\cdot x'}{jA_j^{\frac{1-\beta}{2}}r_j^{1+\beta}}\le 
C\frac{\tilde{R}^{1+\beta}}{j^{\frac{1-\beta}{2}}}\to 0
$$
as $j\to\infty$.
For the second part we use  (\ref{onesidedfrombelowinprop2}): 
$$
\inf_{B_{\tilde{R}}'}h_j(x',0)\ge \inf_{B'_{\sqrt{j}r_j\tilde{R}}}\frac{h_j^+(x',0)}{jA_j^{\frac{1-\beta}{2}}r_j^{1+\beta}}\ge 
- C\frac{\tilde{R}^{1+\beta}}{j^{\frac{1-\beta}{2}}}\to 0\textrm{ as } j\to\infty.
$$

In order to prove the proposition we want to  show that $h_j\to h_0\equiv  0$,
which would clearly contradict uniform convergence and the fact that by definition $\sup_{B_1} |h_0|= \lim_j\sup_{B_1}|h_j|=1$. 

To show the convergence to $0$ we need to control the right-hand side
in equation (\ref{eq:HJforhj}). In particular, we are going to prove that $T^j_1,..., T^j_5\to 0$ as $j\to \infty$.
We formulate this as a new claim.

\vspace{3mm}
\noindent
\textbf{Claim 3:} {\sl We have $T^j_i\to 0$ locally uniformly as $j\to\infty$, for $i=1,...,5$, 
where $T^j_i$ has been defined in equation (\ref{eq:HJforhj}).}

\noindent
\textsl{Proof of Claim 3:}
Using that $|b(x)|\le A|x|^\alpha$ and that $\frac{\alpha}{2 + \alpha}  \le \beta \le \alpha/2$, we see that
$$
\sup_{B_1}|T^j_1|=\sup_{B_1}\frac{2|b|}{\sqrt{j}r_j^{\beta}A_j^{\frac{1-\beta}{2}}}\le
\frac{2A_j\big( A_j^{-\frac{(1-\beta)}{2}}r_j^{1-\beta}\big)^\alpha}
{\sqrt{j}A_j^{\frac{1-\beta}{2}}r_j^{\beta}}=
\frac{2A_j^{\frac{1-\alpha+\beta+\alpha\beta}{2}}r_j^{\alpha-\beta(1+\alpha)}}{\sqrt{j}}\to 0
$$
as $j\to\infty$,
which implies the claim for $T^j_1$.

We also have, by assumption, that $b(0,x_n)\cdot e_n + g_j'(x_n) = 0$,
and thus
\begin{equation}\label{estimatingT2}
\sup_{B_1}|b(x)\cdot e_n-g_j'|\le \sup_{|x_n|\le A_J^{-\frac{1-\beta}{2}}r_j^{1-\beta}}\osc_{x'\in B_{\sqrt{j}r_j}}
\big( b\cdot e_n\big)\le A_jj^{\alpha/2}r_j^\alpha.
\end{equation}
Using Lemma \ref{lemmaA} we see that
\begin{equation}\label{usinglemmaa}
\sup_{B_1}|\partial_n h_j|\le \sup_{(x'/(\sqrt{j}r_j),A_j^{\frac{1-\beta}{2}}x_n/(r_j^{1-\beta}))\in B_1}
\frac{|\partial_n h_j^+ +g_j'|}{jA_j^{1-\beta}r_j^{2\beta}}\le
\frac{c A_j^{-\frac{(1-\beta)^2}{2}}r_j^{-\beta-\beta^2}}{j}.
\end{equation}
From (\ref{estimatingT2}),  (\ref{usinglemmaa}) and $\beta\le \alpha/2\le 1/2$ we see that
$$
\sup_{B_1}|T^j_2|\le C\frac{A_j^{\frac{1+\beta^2}{2}}r_j^{\alpha-\beta-\beta^2}}{j^{1-\frac{\alpha}{2}}} \to 0
$$
as $j\to\infty$.
Next, we estimate in $B_1$
\begin{equation}\label{t3est}
|T^j_3|\le\frac{2|b|^2}{jA_j^{1-\beta}r_j^{2\beta}}\le
\frac{2A_j^{2-(1-\beta)\alpha} r_j^{2\alpha(1-\beta)}}{jA_j^{1-\beta}r_j^{2\beta}}=
\frac{2A_j^{1+\beta-(1-\beta)\alpha}r_j^{2(\alpha-\beta(1+\alpha))}}{j}\to 0
\end{equation}
as $j\to\infty$.
The term $T^j_4$ satisfies
$$
|T^j_4|\le \frac{|b\cdot e_n+ g'_j|}{jA_j^{1-\beta}r_j^{2\beta}}\le
\sup_{|x_n|\le r_j^{1-\beta}}\osc_{x'\in B_{\sqrt{j}r_j}}
\Big(  \frac{b\cdot e_n}{jA_j^{1-\beta}r_j^{2\beta}}  \Big) \le
\frac{A_j^\beta r_j^{\alpha-2\beta}}{j^{1-\alpha/2}}\to 0
$$
as $j\to\infty$.

So far we have shown that $T^j_1,\dots, T^j_4\to 0$ as $j\to \infty$. 

Next, we derive estimates for $T^j_5$. First we use (\ref{relationrjandj}) to arrive at  
$$
\{x;\; (x'/(\sqrt{j}r_j),A_j^{(1-\beta)/2}x_n/r_j^{1-\beta})\in 
B_1 \}\subset B_{A_j^{-(1-\beta)/2}r_j^{1-\beta}}.
$$
Together with Lemma \ref{lemmaA} this implies that
$$
|\tilde\nabla h_j^+|\le CA_j^{\frac{1-\beta}{2}}A_j^{-\frac{1-\beta}{2}\beta}r_j^{(1-\beta)\beta} \quad \textrm{in } 
\{x;\; (x'/(\sqrt{j}r_j),A_j^{\frac{1-\beta}{2}}x_n/r_j^{1-\beta})\in B_1 \},
$$
that is
\begin{equation}\label{estimatesonhjprime}
|\tilde\nabla h_j|\le C\frac{ A_j^{-\frac{1-\beta}{2}\beta}r_j^{-\beta^2}}{\sqrt{j}},
\end{equation}

Now, since $h_j$ is a Lipschitz solution to the Hamilton-Jacobi equation,
the super-differential consists almost everywhere of a unique element $p$,
and almost everywhere, $\nabla h_j =p$ solves equation (\ref{eq:HJforhj}). That is, at
almost every point,
$$
p_n^2-\frac{2}{jA_j^{1-\beta}r_j^{2\beta}}p_n+\frac{\delta}{jA_j^{1-\beta}r_j^{2\beta}}=0,
$$
where
$$
\delta=|\tilde\nabla h_j|^2-\frac{2b}{\sqrt{j}A_j^{\frac{1-\beta}{2}}r_j^{\beta}}\cdot\tilde\nabla h_j
-2\partial_n h_j ((b(x)\cdot e_n+ g_j')
$$
$$
+\frac{|b + e_n g_j'|^2}{jA_j^{1-\beta}r_j^{2\beta}}+2\frac{e_n\cdot b + g_j'}{jA_j^{1-\beta}r_j^{2\beta}}.
$$
This means that $p_n$ can only take the two values (depending on $\delta$),
\begin{equation}\label{livetarskit1}
p_n=
\left\{
\begin{array}{l}
\frac{1}{jA_j^{1-\beta}r_j^{2\beta}}+\sqrt{\frac{1}{j^2A_j^{2-2\beta}r_j^{4\beta}}-\frac{\delta}{jA_j^{1-\beta}r_j^{2\beta}}} \\
\frac{1}{jA_j^{1-\beta}r_j^{2\beta}}-\sqrt{\frac{1}{j^2A_j^{2-2\beta}r_j^{4\beta}}-\frac{\delta}{jA_j^{1-\beta}r_j^{2\beta}}}.
\end{array}
\right.
\end{equation}
By (\ref{usinglemmaa}) as well as the above estimates on
$|\tilde\nabla h_j|$ (c.f. equation (\ref{estimatesonhjprime}))  we see that
\begin{equation}\label{livetarskit2}
|\delta|\le C\Bigg[ \frac{A_j^{-(1-\beta)\beta}r_j^{-2\beta^2}}{j}+
\frac{ A_j^{\frac{1-\alpha+\alpha\beta+\beta^2}{2}}r_j^{\alpha-\beta(1+\alpha+\beta)}}{j}
\end{equation}
$$
+\frac{A_j^{\frac{1+2\beta-\beta^2}{2}}r_j^{\alpha-\beta-\beta^2}}{j^{1-\frac{\alpha}{2}}} +
\frac{A_j^{1+\beta-(1-\beta)\alpha}r_j^{2(\alpha-\beta(1+\alpha))}}{j}+\frac{A_j^\beta r_j^{\alpha-2\beta}}{j^{1-\frac{\alpha}{2}}}
\Bigg].
$$
Next, using (\ref{relationrjandj}) and (\ref{t3est}) we see that when $j$ is large
and $r_j$ is small,
\begin{equation}\label{livetarskit3}
\frac{A_j^{\frac{1+2\beta-\beta^2}{2}}r_j^{\alpha-\beta-\beta^2}}{j^{1-\frac{\alpha}{2}}}\le
\Big(A_j^{\frac{1-\beta^2}{2}}r_j^{\beta-\beta^2} \Big) 
\frac{A_j^{\beta}r_j^{\alpha-2\beta}}{j^{1-\frac{\alpha}{2}}}< < 
\frac{A_j^{\beta}r_j^{\alpha-2\beta}}{j^{1-\frac{\alpha}{2}}}.
\end{equation}
Similarly, using  (\ref{relationrjandj}) which implies that
$$
r_j^{\beta-\alpha\beta-\beta^2}\le j^{\frac{\alpha+\beta-1}{1+\beta}}A_j^{-\frac{(1-\beta)(1-\alpha-\beta)}{2}},
$$
we may conclude that
\begin{equation}\label{livetarskit224}
\frac{A_j^{\frac{1-\alpha+\alpha\beta+\beta^2}{2}}r_j^{\alpha-\beta(1+\alpha)-\beta^2}}{j}=
\Big(\frac{A_j^{\frac{1-\alpha+\alpha\beta+\beta^2-2\beta}{2}}}{j^{\alpha/2}}r_j^{\beta-\alpha\beta-\beta^2} \Big)\frac{A_j^\beta r_j^{\alpha-2\beta}}{j^{1-\alpha/2}}
\end{equation}
$$
\le j^{\frac{\alpha+\beta-1}{1+\beta}-\frac{\alpha}{2}}\frac{A_j^\beta r_j^{\alpha-2\beta}}{j^{1-\alpha/2}},
$$
Noticing that the term $j^{\frac{\alpha+\beta-1}{1+\beta}-\frac{\alpha}{2}}$ in (\ref{livetarskit224})
goes to zero as $j\to \infty$,
we deduce that
\begin{equation}\label{livetarskit4}
\frac{A_j^{\frac{1-\alpha+\beta+2\alpha\beta+\beta^2}{2}}r_j^{\alpha-\beta(1+\alpha)-\beta^2}}{j}< < \frac{A_j^\beta r_j^{\alpha-2\beta}}{j^{1-\alpha/2}}\textrm{ as }j\to \infty.
\end{equation}
%% OLD OLD OLD
%%\begin{equation}
%%=\Big( j^{-\frac{1-\beta}{2}}
%%A_j^{\frac{1+\beta+\alpha\beta-\alpha}{2}}\Big)
%%\frac{A_j^\beta r_j^{\alpha-2\beta}}{j^{1-\frac{\alpha}{2}}}< < 
%%\frac{A_j^\beta r_j^{\alpha-2\beta}}{j^{1-\frac{\alpha}{2}}}.
%%\end{equation}
Moreover, as $2\beta\le \alpha\le 1$, we have 
\begin{equation}\label{livetarskit5}
\frac{A_j^{1+\beta-(1-\beta)\alpha}r_j^{2(\alpha-\beta(1+\alpha))}}{j}< < \frac{A_j^{-(1-\beta)\beta}r_j^{-2\beta^2}}{j}
\end{equation}
when $r_j$ is small. Using (\ref{livetarskit2})-(\ref{livetarskit5}), we estimate
$$
|\delta|\le C\left[ \frac{A_j^{-(1-\beta)\beta}r_j^{-2\beta^2}}{j}+
\frac{A_j^\beta r_j^{\alpha-2\beta}}{j^{1-\frac{\alpha}{2}}}
\right].
$$
Using (\ref{relationrjandj}) again along with  $\alpha/(2+\alpha)\le \beta\le \alpha/2\le 1/2$, we see that
\begin{equation}\label{livetarskit6}
\frac{A_j^{-(1-\beta)\beta} r_j^{-2\beta^2}}{j} = \Big( \frac{1}{j^{2\frac{\beta-\beta^2}{\beta+\beta^2}}}\Big)\frac{A_j^{-(1-\beta)}r_j^{-2\beta}}{j} < <
\frac{A_j^{-(1-\beta)}r_j^{-2\beta}}{j}
\end{equation}
and that
\begin{equation}\label{livetarskit7}
\frac{A_j^\beta r_j^{\alpha-2\beta}}{j^{1-\frac{\alpha}{2}}}\le \Big( \frac{A_j^{\frac{2\beta-\alpha+\alpha\beta}{2\beta}}}{j^{\frac{-\alpha(\beta+\beta^2)+2\alpha}{2(\beta+\beta^2)}}}\Big)\frac{A_j^{-(1-\beta)}r_j^{-2\beta}}{j} < < \frac{A_j^{-(1-\beta)}r_j^{-2\beta}}{j}
\end{equation}
when $j$ is large.
From (\ref{livetarskit1}), (\ref{livetarskit6}) and (\ref{livetarskit7}) we see that we can make a 
Taylor expansion in equation (\ref{livetarskit1}) and deduce that $p_n$ is of order
\begin{equation}\label{twoforpn}
|p_n|\approx
\left\{
\begin{array}{l}
\frac{2}{jA_j^{1-\beta}r_j^{2\beta}} \\
\frac{|\delta|}{2}\le \frac{CA_j^{-(1-\beta)\beta}r_j^{-2\beta^2}}{j}+
\frac{CA_j^\beta r_j^{\alpha-2\beta}}{j^{1-\frac{\alpha}{2}}}.
\end{array}
\right.
\end{equation}
Now remember that from (\ref{usinglemmaa}), we have the additional information that
$$
|\partial_n h_j|\le \frac{c A_j^{-\frac{1-2\beta+\beta^2}{2}} r_j^{-\beta-\beta^2}}{j}.
$$
As $\beta\le 1/2$, we conclude that 
for small $r_j$, the first line in (\ref{twoforpn}) and not the second
must hold.
Thus
$$
|p_n|= \frac{|\delta|}{2}+
\textrm{lower order terms.}
$$
It follows that almost everywhere,
$$
T_5=jA_j^{1-\beta}r_j^{2\beta}(p_n)^2\le C\left( \frac{A_j^{1-3\beta+2\beta^2}r_j^{2\beta-4\beta^2}}{j}+
\frac{A_j^{1+\beta}r_j^{2(\alpha-\beta)}}{j^{1-\alpha}}
\right).
$$
Since $2\beta\le \alpha\le 1$, the right-hand
goes to zero as $j\to \infty$, and the claim follows.

The proof of the proposition is now easy to finish. We know from Claim 3 and Lemma \ref{ComparisonParabolic} that
$h_j\to h_0$ locally uniformly in $\R^n$, 
 where $h_0$ is a viscosity solution to
\begin{equation}\label{ujk}
|\nabla ' h_0|^2-2\partial_n h_0 = 0.
\end{equation}
Moreover $h_j\to 0$ on $\{x_n=0\}$ by Claim 2. Using uniqueness for
parabolic Hamilton-Jacobi equations (\cite{L}) we deduce that $h_0=0$
for $x_n\le 0$ and $h_0\ge 0$ in $\{x_n\ge 0\}$.
By Claim 2 we also know that $h_0$ is concave in the $x'$ directions, which
together with the bound from below implies that $h_0$ is constant for each $x_n\ge 0$: if $h_0(x',x_n)$ is not constant for some
$x_n=t> 0$, then $h_0(x',t)\le h_0(0,t)+p\cdot x'$ for some non-zero $p\in \mathbb{R}^{n-1}$. It follows that there is a point $x_0'$ such that
$h_0(x_0',t)<0$. But from (\ref{ujk}) it follows that $\partial_n h_0\ge 0$, a contradiction. Therefore $h_0(x)=h_0(x_n)$.
Using (\ref{ujk}) again we see that $\partial_n h_0=0$ in $\{x_n\ge 0\}$, and thus $h_0=0$. This contradicts $\sup_{B_1}|h_0|=1$, and the Proposition follows. \qed

%%%%%%%%%%%%%%%%%%%%%%%%%%%%%%%%%%%%%%

%%%%%%%%%%%%%%%%%%%%%%%%%%%%%%%%%%%%%%%%%%

\section{Improved Regularity in Degenerate Directions}

In this section we improve the regularity in the directions orthogonal to the
non-degenerate direction. The idea is to use the scaling from Proposition \ref{nondeg}
to get a parabolic equation and deduce better regularity in the $x'$ directions
from the regularity theory for parabolic Hamilton-Jacobi equations.

\begin{prop}\label{orthogonal} Let $h^+$ be a viscosity solution to
$|\nabla h^+-a|^2=1$ in $B_R$ and assume that $h^+$
satisfies the one-sided $C^{1,\beta}$ estimate from above
\begin{equation}
\sup_{B_r(0)}\big( h^+(x) - h^+(0)-p\cdot (x-0)\big)\le  C A^{1-\beta\over 2}r^{1+\beta}.
\end{equation}
Furthermore let
\begin{equation}
h^+(x',x_n)\le h^+(0,x_n)+p\cdot x' + CA^{\frac{1-\beta}{2}} |x'|^{1+\beta}
\end{equation}
for each fixed $x_n$, and assume that $h^+$ satisfies the one-sided $C^{1,\beta}$-estimate from below at the origin
\begin{equation}
h^+(x)\ge -CA^{\frac{1-\beta}{2}}|x|^{1+\beta}.
\end{equation}
Also assume that
\begin{equation}
h^+(x',0)\ge h^+(0)+p\cdot x' - CA^{\frac{1-\beta}{2}} |x'|^{1+\beta}.
\end{equation}
In the just mentioned estimates
we assume moreover that the constant $A\le \tilde{C}$ (where $\tilde{C}$ is a fixed constant) 
controls the 
seminorm of $a(x)$ by $[a]_{C^\alpha}\le A$. 
Finally we assume that $a(x)=e_n+b(x)$ where $b(0)=0$ and that $\alpha/(2+\alpha)\le \beta\le \alpha/2$.

Then for every every $\delta\le\min(A^{1/2}R^{\frac{1}{1-\beta}},R)$ and every 
$x_n\in (0,A^{-\frac{1-\beta}{2}}\delta^{1-\beta})$ we have for some $p\in \mathbb{R}^n$
$$
\sup_{B_\delta'(0,x_n)}\big( h^+(x)-h^+(0,x_n)-p\cdot x'\big)\le  Q
$$
where
$$
Q=\left\{\begin{array}{ll}
CA^{\frac{1}{2}-\frac{\alpha\beta}{2(2+\alpha-2\beta)}}\delta^{1+\frac{\alpha}{2+\alpha-2\beta}} &
\textrm{if } \delta\le \min\big( A^{\frac{6\beta+\alpha\beta-2-\alpha}{2(4\beta-\alpha)}},
A^{-\frac{\beta}{\alpha-2\beta}},A^{\frac{\beta}{2}}R^{\frac{2+\alpha-2\beta}{2}}\big), \\
CA^{\frac{1-\beta}{2}}R^{\beta-1}\delta^2 & \textrm{if } A^{\frac{\beta}{2}}R^{\frac{2+\alpha-2\beta}{2}}\le \delta\le
A^{\frac{6\beta+\alpha\beta-2-\alpha}{2(4\beta-\alpha)}}, \\
CA^{\frac{1+\beta}{2}}\delta^{1+\alpha-\beta} & \textrm{if } A^{-\frac{\beta}{\alpha-2\beta}}\le \delta \le
A^{\frac{6\beta+\alpha\beta-2-\alpha}{2(4\beta-\alpha)}}, \\
CA^{\frac{1-\beta}{1+\beta}}\delta^{1+\frac{2\beta}{1+\beta}} & \textrm{if }A^{\frac{6\beta+\alpha\beta-2-\alpha}{2(4\beta-\alpha)}}\le \delta\le \min\big( A^{-\frac{1-\beta}{2\beta}},A^{\frac{1-\beta}{2}}R^{1+\beta}\big),\\
CA^{\frac{3(1-\beta)}{2}}\delta^{1+3\beta} & \textrm{if } \min\bigg[\max\big( A^{-\frac{1-\beta}{2\beta}},A^{\frac{6\beta+\alpha\beta-2-\alpha}{2(4\beta-\alpha)}}\big),\\ &\quad \quad\max\big(A^{\frac{1-\beta}{2}}R^{1+\beta},A^{\frac{6\beta+\alpha\beta-2-\alpha}{2(4\beta-\alpha)}}\big)\bigg]\le \delta
\end{array}
\right.
$$
if $R=\min(A^{1/2}R^{\frac{1}{1-\beta}},R)$, and
$$
Q=\left\{ \begin{array}{ll}
CA^{\frac{1}{2}-\frac{\alpha\beta}{2(2+\alpha-2\beta)}}\delta^{1+\frac{\alpha}{2+\alpha-2\beta}} &
\textrm{if } \delta\le \min\big(A^{\frac{6\beta+\alpha\beta-2-\alpha}{2(4\beta-\alpha)}},A^{-\frac{\beta}{\alpha-2\beta}},A^{\frac{2+\alpha}{4}}R^{\frac{2+\alpha-2\beta}{2(1-\beta)}}\big), \\
CR^{-1}\delta^2 & \textrm{if } A^{\frac{2+\alpha}{4}}R^{\frac{2+\alpha-2\beta}{2(1-\beta)}}\le \delta\le A^{\frac{6\beta+\alpha\beta-2-\alpha}{2(4\beta-\alpha)}}, \\
CA^{\frac{1+\beta}{2}}\delta^{1+\alpha-\beta} & \textrm{if } A^{-\frac{\beta}{\alpha-2\beta}}\le \delta \le A^{\frac{6\beta+\alpha\beta-2-\alpha}{2(4\beta-\alpha)}}, \\
CA^{\frac{1-\beta}{1+\beta}}\delta^{1+\frac{2\beta}{1+\beta}} & \textrm{if } A^{\frac{6\beta+\alpha\beta-2-\alpha}{2(4\beta-\alpha)}}\le \delta\le AR^{\frac{1+\beta}{1-\beta}}, \\
CR^{-1}\delta^2 & \textrm{if } \max\big( AR^{\frac{1+\beta}{1-\beta}},A^{\frac{6\beta+\alpha\beta-2-\alpha}{2(4\beta-\alpha)}}\big)\le \delta 
\end{array}\right.
$$
if $A^{1/2}R^{\frac{1}{1-\beta}}=\min(A^{1/2}R^{\frac{1}{1-\beta}},R)$.
\end{prop}

\textsl{Proof:} From Proposition \ref{nondeg} we obtain that
$$
h(x)=
\frac{h^+(rx',A^{-(1-\beta)/2}r^{1-\beta}x_n)+g(A^{-(1-\beta)/2}r^{1-\beta}x_n)}{A^{(1-\beta)/2}r^{1+\beta}}
$$
is bounded from above. Moreover we know from the proof of 
Proposition \ref{nondeg} that $h$ solves
(with slightly reordered terms as compared to equation (\ref{eq:HJforhj}))
\begin{align}\label{rescaledequation}
\big|\tilde\nabla h-\frac{b'(0,A^{-(1-\beta)/2}r^{1-\beta}x_n)}{A^{(1-\beta)/2}r^\beta} \big|^2 - 2\partial_n h
\\\nonumber =
-A^{1-\beta}r^{2\beta}|\partial_n h|^2+
2\frac{b(rx',A^{-(1-\beta)/2}r^{1-\beta}x_n)-b(0,A^{-(1-\beta)/2}r^{1-\beta}x_n)}{A^{(1-\beta)/2}r^\beta}
\cdot \tilde\nabla h
\\\nonumber 
-2\frac{b(rx',A^{-(1-\beta)/2}r^{1-\beta}x_n)\cdot e_n+g'}{A^{1-\beta}r^{2\beta}}+
2\big( b(rx',A^{-(1-\beta)/2}r^{1-\beta}x_n)\cdot e_n+g' \big)\partial_n h
\\\nonumber -
\frac{|b(rx',A^{-(1-\beta)/2}r^{1-\beta}x_n)+g'e_n|^2}{A^{1-\beta}r^{2\beta}}-
\frac{|b'(0,A^{-(1-\beta)/2}r^{1-\beta}x_n)|^2}{A^{1-\beta}r^{2\beta}}
\\\nonumber =
-I-II-III-IV-V-VI,
\end{align}
where we have used the notation $b'=(b_1,b_2,...,b_{n-1},0)$ in the first line of the equation.

Let us first show that the right-hand is bounded:
Using exactly the same argument as in the estimate of $T_5$ in Claim 3 of Proposition \ref{nondeg}
(consider the
same argument with $j:= 1$), we see that 
\begin{align}\label{est-ii}
    \sup_{B_1}|I|&\le C(A^{1-3\beta + 2\beta^2}r^{2\beta(1-2\beta)} + A^{1+ \beta}r^{2(\alpha - \beta)}) \leq C .
\intertext{Next,}
\label{est-iii}
\sup_{B_1}|II| &\le \sup_{B_1}
\frac{|b(rx',A^{-(1-\beta)/2}r^{1-\beta}x_n)-b(0,A^{-(1-\beta)/2}r^{1-\beta}x_n)|}{A^{(1-\beta)/2}r^\beta}
|\tilde\nabla h|\\ \nonumber
&\le 2 A^{\frac{1+\beta}{2}}r^{\alpha-\beta}\sup_{B_1}|\tilde\nabla h(x',x_n)|.
\intertext{Similarly we may estimate}
\label{est-iv}
    \sup_{B_1}|III| &\le \sup_{|x_n|\le A^{\frac{-1+\beta}{2}}r^{1-\beta}}\osc_{B'_r}
    \Big( \frac{b\cdot e_n}{A^{1-\beta}r^{2\beta}}\Big)\le 4A^{\beta}r^{\alpha-2\beta} \le C,
\\
\label{est-v}
    \sup_{B_1}|IV|& \le 2 A^{(1+\beta)/2}r^{\alpha-\beta}\osc_{B_1}\sqrt{I} \le C,
\\
\label{est-i}
    \sup_{B_1} |V| &\leq CA^{(1-\alpha) + (1+ \alpha)\beta}r^{2\alpha - 2\beta(1+\alpha)} \leq C,\\
\label{est-vi}
    \sup_{B_1} |VI| &\leq C A^{(1-\beta)(1-\alpha)} r^{2(\alpha-\beta-\alpha\beta)} \leq C.
\end{align}

From (\ref{est-ii}), (\ref{est-iii}), (\ref{est-iv}), (\ref{est-v}), (\ref{est-i}) and(\ref{est-vi}) as well as the parabolic comparison principle Lemma \ref{ComparisonParabolic}, we obtain that
$|\nabla h|\le C$ in $Q_1=\{ x;\; |x'|\le 1,\; 0< x_n <1\}$.

Next we are going to 
estimate the oscillation of each term of the right-hand side in
equation (\ref{rescaledequation})
in the $x'$ variable:
First we notice that the oscillation of $VI$ in the $x'$
directions is identically zero. 

Next, using the gradient bound of $h$, we obtain that for every $0< x_n <1$,
\begin{align}
\osc_{B_1'}I\le CA^{1-\beta}r^{2\beta}\label{osc1},\\
\osc_{B_1'}II\le  2A^{(1+\beta)/2}r^{\alpha-\beta}\label{osc2},\\
 \osc_{B_1'}III \le
4A^{\beta}r^{\alpha-2\beta}\label{osc3},\\
\osc_{B_1'}IV\le CAr^{\alpha}\label{osc4},\\
\osc_{B_1'} V\le 2A^{1+\beta}r^{2(\alpha-\beta)}.\label{osc5}
\end{align}

So the oscillation in $x'$ of the right-hand side in equation (\ref{rescaledequation})
is estimated by $q=C\max\big(A^{1-\beta}r^{2\beta},A^{\beta}r^{\alpha-2\beta}\big)$, 
where we have used that in view of $2\beta\le \alpha \le 1$, $\osc_{B_1'}I$ or $\osc_{B_1'}III$ is dominating when $r$ is small.

Now let $w$ be the solution to
$$
\begin{array}{ll}
|\tilde\nabla w|^2-2w_n=q+f(x_n) & \textrm{in } Q_1, \\
w=h & \textrm{on } \partial{Q_1}\setminus \{x_n=0\},
\end{array}
$$
where $f(x_n)$ is chosen
so that the right-hand side of the above equation equals the right-hand side of (\ref{rescaledequation})
at $|x'|=0$.
By the parabolic comparison principle Lemma \ref{ComparisonParabolic} and the one-sided estimate Lemma \ref{OneSidedEstimates}
for viscosity solutions we have
for $p$ in the super-differential of $w$ at $x_0$ that
$$
h(x)\le w(x) \le w(x_0)+p\cdot (x'-x'_0)+C|x'-x'_0|^2
$$
$$\le
h(x_0)+p\cdot (x'-x'_0)+C|x'-x'_0|^2+Cq.
$$
Rescaling back to $h^+$ we see that this implies that (for some $p$ in the super-differential of $h^+$)
$$
h^+(x',x_n)\le h^+(0,x_n)+p\cdot x'+C
\big( A^{(1-\beta)/2}r^{1+\beta}q+A^{(1-\beta)/2}r^{\beta-1}|x'|^2 \big)
$$
for $0<x_n<A^{-{{1-\beta}\over 2}}r^{1-\beta}$ and $|x'|\le r$.

Using this estimate in an optimal way will yield the proposition. Rearranging terms and taking the
supremum in $B_\delta '$ on both sides leads to
\begin{equation}\label{xprimeestimatefromabove}
\sup_{B_\delta'}\big( h^+(x',x_n)- h^+(0,x_n)-p\cdot x' \big)\le C
\big( A^{(1-\beta)/2}r^{1+\beta}q+A^{(1-\beta)/2}r^{\beta-1}\delta^2 \big)
\end{equation}
provided that $\delta\in (0,r)$ and $0<x_n<A^{-{{1-\beta}\over 2}}\delta^{1-\beta}$.

For fixed $R$ and fixed $\delta \in (0,R)$ we 
want to find an $r$ optimizing  
this estimate. Since we have a constraint that $r\le \min(A^{1/2}R^{\frac{1}{1-\beta}},R)$ 
it is natural to divide the
proof into two cases.

\textbf{Case 1, $R=\inf(A^{1/2}R^{\frac{1}{1-\beta}},R)$:}
That is, when
$$
R\ge A^{-\frac{1-\beta}{2\beta}}.
$$
When $q=CA^{\beta}r^{\alpha-2\beta}$, or equivalently when 
\begin{equation}\label{thu1}
r\le A^{-\frac{1-2\beta}{4\beta-\alpha}},
\end{equation} 
then we would want to choose 
\begin{equation}\label{thu2}
r=A^{-\frac{\beta}{2+\alpha-2\beta}}\delta^{\frac{2}{2+\alpha-2\beta}}.
\end{equation}
From (\ref{thu1}) we infer then that
\begin{equation}\label{thu3}
A^{-\frac{\beta}{2+\alpha-2\beta}}\delta^{\frac{2}{2+\alpha-2\beta}}\le A^{-\frac{1-2\beta}{4\beta-\alpha}},
\end{equation}
or equivalently that
\begin{equation}\label{thu4}
\delta\le A^{\frac{6\beta+\alpha\beta-2-\alpha}{2(4\beta-\alpha)}}.
\end{equation}
Moreover we need $r\ge \delta$ and $r\le R$ in order to use (\ref{xprimeestimatefromabove}), that is
\begin{equation}\label{thu5}
\delta\le A^{-\frac{\beta}{\alpha-2\beta}}
\end{equation}
and
\begin{equation}\label{thu6}
\delta\le A^{\frac{\beta}{2}}R^{\frac{2+\alpha-2\beta}{2}}.
\end{equation}
When (\ref{thu1}), (\ref{thu5}) and (\ref{thu6}) hold, then we can choose $r$ according to 
(\ref{thu2}) and deduce that
$$
\sup_{B_{\delta}'}\big( h^+(x',x_n)-h^+(0,x_n)-p\cdot x' \big)\le
CA^{\frac{1}{2}-\frac{\alpha\beta}{2(2+\alpha-2\beta)}}\delta^{1+\frac{\alpha}{2+\alpha-2\beta}}.
$$

It still remains to consider the cases when at least one of (\ref{thu1}), (\ref{thu5}) or 
(\ref{thu6}) does not hold.

If (\ref{thu1}) holds but not (\ref{thu6}) then we choose 
$r=R$ and we deduce that ---using $q=A^\beta r^{\alpha-2\beta}$ which is equivalent to (\ref{thu3})---
\begin{equation}\label{thu7}
\sup_{B_{\delta}'}\big( h^+(x',x_n)-h^+(0,x_n)-p\cdot x' \big)\le
C\big( A^{\frac{1-\beta}{2}}R^{1+\beta}q+A^{\frac{1-\beta}{2}}R^{\beta-1}\delta^2\big)
\end{equation}
$$
=C\big(A^{\frac{1+\beta}{2}}R^{1+\alpha-\beta}+A^{\frac{1-\beta}{2}}R^{\beta-1}\delta^2 \big).
$$
In 
order to simplify (\ref{thu7})
we use that $\delta > A^{\frac{\beta}{2}}R^{\frac{2+\alpha-2\beta}{2}}$ and thus
$$
A^{\frac{1-\beta}{2}}R^{\beta-1}\delta^2\ge A^{\frac{1+\beta}{2}}R^{1+\alpha-\beta}.
$$
This way we may simplify (\ref{thu7}) to
\begin{equation}\label{thu8}
\sup_{B_{\delta}'}\big( h^+(x',x_n)-h^+(0,x_n)-p\cdot x' \big)\le CA^{\frac{1-\beta}{2}}R^{\beta-1}\delta^2.
\end{equation}

If (\ref{thu1}) holds but (\ref{thu5}) does not, then we choose $r=\delta$ and deduce
from (\ref{xprimeestimatefromabove}) that
\begin{equation}\label{thu9}
\sup_{B_{\delta}'}\big( h^+(x',x_n)-h^+(0,x_n)-p\cdot x' \big)\le
C\big( A^{\frac{1+\beta}{2}}\delta^{1+\alpha-\beta}+A^{\frac{1-\beta}{2}}\delta^{1+\beta}\big).
\end{equation}
The information that (\ref{thu5}) does not hold implies
$$
\delta\ge A^{-\frac{\beta}{\alpha-\beta}},
$$
so we may estimate
$$
A^{\frac{1+\beta}{2}}\delta^{1+\alpha-\beta}=\Big( A^\beta\delta^{\alpha-2\beta}\Big)A^{\frac{1-\beta}{2}}\delta^{1+\beta}\ge A^{\frac{1-\beta}{2}}\delta^{1+\beta}.
$$
We can thus simplify (\ref{thu9}) to
\begin{equation}\label{thu10}
\sup_{B_{\delta}'}\big( h^+(x',x_n)-h^+(0,x_n)-p\cdot x' \big)\le
C A^{\frac{1+\beta}{2}}\delta^{1+\alpha-\beta}.
\end{equation}

In case (\ref{thu1}) does not hold ---that is when $q=CA^{1-\beta}r^{2\beta}$--- our
first choice of $r$ is
$r=A^{-\frac{1-\beta}{2(1+\beta)}}\delta^{\frac{1}{1+\beta}}$; notice that by the condition
$r\le R$,
this is only possible for 
\begin{equation}\label{thu11}
\delta\le A^{\frac{1-\beta}{2}}R^{1+\beta}.
\end{equation}
Moreover we must have that $r\ge \delta$ which in this case becomes 
\begin{equation}\label{thu12}
\delta\le A^{-\frac{1-\beta}{2\beta}}.
\end{equation}
So in case (\ref{thu1}) does not hold but (\ref{thu11}) and (\ref{thu12}) do we obtain
from (\ref{xprimeestimatefromabove}) that
\begin{equation}\label{thu13}
\sup_{B_{\delta}'}\big( h^+(x',x_n)-h^+(0,x_n)-p\cdot x' \big)\le
C A^{\frac{1-\beta}{1+\beta}}\delta^{1+\frac{2\beta}{1+\beta}}.
\end{equation}

If neither (\ref{thu1}) nor (\ref{thu12}) hold then we chose $r=2\delta$
and deduce that
\begin{equation}\label{thu14}
\sup_{B_{\delta}'}\big( h^+(x',x_n)-h^+(0,x_n)-p\cdot x' \big)\le
C \big( A^{\frac{3(1-\beta)}{2}}\delta^{1+3\beta}+A^{\frac{1-\beta}{2}}\delta^{1+\beta}\big).
\end{equation}
In order to simplify (\ref{thu14}) we notice that as (\ref{thu12}) does not hold, we have 
$\delta\ge A^{-\frac{1-\beta}{2\beta}}$ and thus
$$
A^{\frac{3(1-\beta)}{2}}\delta^{1+3\beta}\ge A^{\frac{1-\beta}{2}}\delta^{1+\beta}.
$$
Therefore we may write (\ref{thu14}) as
\begin{equation}\label{thu15}
\sup_{B_{\delta}'}\big( h^+(x',x_n)-h^+(0,x_n)-p\cdot x' \big)\le C A^{\frac{3(1-\beta)}{2}}\delta^{1+3\beta}.
\end{equation}

Finally, if neither (\ref{thu1}) nor (\ref{thu11})
hold,
we use $r=R$ in equation (\ref{xprimeestimatefromabove}) 
and $\delta\ge A^{\frac{1-\beta}{2}}R^{1+\beta}$
to obtain
$$
\sup_{{B'_\delta}}\big( h^+ -h^+(0,x_n)-p\cdot x'\big)\le C\big(
A^{\frac{1+\beta}{2}}R^{1+\alpha-\beta}+A^{\frac{1-\beta}{2}}R^{\beta-1}\delta^2\big)\le
CA^{\frac{3(1-\beta)}{2}}\delta^{1+3\beta}.
$$
This is
the final estimate on $Q$ in Case 1.

\textbf{Case 2, $A^{1/2}R^{\frac{1}{1-\beta}}=\inf(A^{1/2}R^{\frac{1}{1-\beta}},R)$:}
The argument is similar as in Case 1. In this case we have $R\le A^{-\frac{1-\beta}{2\beta}}$.
Again, if $q=CA^\beta r^{\alpha-2\beta}$ then the best balance in (\ref{xprimeestimatefromabove}) 
is obtained by the choice $r=A^{-\frac{\beta}{2+\alpha-2\beta}}\delta^{\frac{2}{2+\alpha-2\beta}}$.
But when $q=CA^\beta r^{\alpha-2\beta}$, then $r\le A^{-\frac{1-2\beta}{4\beta-\alpha}}$ which is 
equivalent to
\begin{equation}\label{thu15}
\delta\le A^{\frac{6\beta+\alpha \beta-2-\alpha}{2(4\beta-\alpha)}}
\end{equation}
when $r=A^{-\frac{\beta}{2+\alpha-2\beta}}\delta^{\frac{2}{2+\alpha-2\beta}}$.
In order to use (\ref{xprimeestimatefromabove}) we need $r\ge \delta$ which with our choice of $r$
is equivalent to
\begin{equation}\label{thu16}
\delta\le A^{-\frac{\beta}{\alpha-2\beta}}.
\end{equation}
Moreover, as $r\le A^{\frac{1}{2}}R^{\frac{1}{1-\beta}}$,  
we must have
\begin{equation}\label{thu17}
\delta\le A^{\frac{2+\alpha}{4}}R^{\frac{2+\alpha-2\beta}{2(1-\beta)}}
\end{equation}
for our choice of $r$.

If (\ref{thu15}), (\ref{thu16}) and (\ref{thu17}) hold, then (\ref{xprimeestimatefromabove})
implies that
\begin{equation}\label{thu18}
\sup_{{B'_\delta}}\big( h^+ -h^+(0,x_n)-p\cdot x'\big)\le C A^{\frac{1}{2}-\frac{\alpha\beta}{2(2+\alpha-2\beta)}}\delta^{1+\frac{\alpha}{2+\alpha-2\beta}}.
\end{equation}

As before we also need to investigate what happens if at least One of the conditions
$q=CA^\beta r^{\alpha-2\beta}$, (\ref{thu16}) or 
(\ref{thu17}) does not hold.

If (\ref{thu15}) holds but (\ref{thu17}) does not, then we choose 
$r=A^{\frac{1}{2}}R^{\frac{1}{1-\beta}}$ and deduce from (\ref{xprimeestimatefromabove})
that
\begin{equation}\label{thu19}
\sup_{{B'_\delta}}\big( h^+ -h^+(0,x_n)-p\cdot x'\big)\le C \big( A^{\frac{1+\beta}{2}}r^{1+\alpha-\beta}
+A^{\frac{1-\beta}{2}}R^{\beta-1}\delta^2\big)
\end{equation}
$$
=C\big( A^{\frac{2+\alpha}{2}}R^{\frac{1+\alpha-\beta}{1-\beta}}+R^{-1}\delta^2\big).
$$
On the other hand, the information that (\ref{thu17})
does not hold implies that
$R^{-1}\delta^2\ge A^{\frac{2+\alpha}{2}}R^{\frac{1+\alpha-\beta}{1-\beta}}$. We may therefore simplify (\ref{thu19}) to 
\begin{equation}\label{thu20}
\sup_{{B'_\delta}}\big( h^+ -h^+(0,x_n)-p\cdot x'\big)\le CR^{-1}\delta^2.
\end{equation}

Next, if (\ref{thu15}) holds but (\ref{thu16}) does not, we choose $r=2\delta$
and deduce that ---using that (\ref{thu16}) does not hold---
\begin{equation}\label{thu21}
\sup_{{B'_\delta}}\big( h^+ -h^+(0,x_n)-p\cdot x'\big)\le C\big( A^{\frac{1+\beta}{2}}\delta^{1+\alpha-\beta}+A^{\frac{1-\beta}{2}}\delta^{1+\beta}\big)
\end{equation}
$$
\le CA^{\frac{1+\beta}{2}}\delta^{1+\alpha-\beta}.
$$

If (\ref{thu15}) does not hold, then $q=CA^{1-\beta}r^{2\beta}$ and the optimal choice of $r$
in (\ref{xprimeestimatefromabove}) is given by
\begin{equation}\label{thu22}
r=A^{-\frac{1-\beta}{2(1+\beta)}}\delta^{\frac{1}{1+\beta}}.
\end{equation}
In order to satisfy $\delta\le r$, we require furthermore that
\begin{equation}\label{thu23}
\delta\le A^{-\frac{1-\beta}{2\beta}}.
\end{equation}
And in order to satisfy $r\le A^{\frac{1}{2}}R^{\frac{1}{1-\beta}}$ we require that
\begin{equation}\label{thu24}
\delta\le AR^{\frac{1+\beta}{1-\beta}}.
\end{equation}
From (\ref{thu23}), (\ref{thu24}), (\ref{thu22}) as well as the
information that (\ref{thu15}) does not hold we infer that

\begin{equation}\label{thu25}
\sup_{{B'_\delta}}\big( h^+ -h^+(0,x_n)-p\cdot x'\big)\le CA^{\frac{1-\beta}{1+\beta}}\delta^{1+\frac{2\beta}{1+\beta}}.
\end{equation}

Notice that (\ref{thu23}) is always satisfied in {\bf Case 2} as
$\delta< R \le A^{-\frac{1-\beta}{2\beta}}$.

We therefore only need to check what happens if neither (\ref{thu24}) nor (\ref{thu15}) hold
in order to finish the proof. In this case we choose $r=A^{\frac{1}{2}}R^{\frac{1}{1-\beta}}$. With
this choice, (\ref{xprimeestimatefromabove}) implies that
---using the fact that (\ref{thu24}) does not hold in the last of the following inequalities---
\begin{equation}\label{thu26}
\sup_{{B'_\delta}}\big( h^+ -h^+(0,x_n)-p\cdot x'\big)\le C \big( A^{\frac{1(1-\beta)}{2}}r^{1+3\beta}+
A^{\frac{1-\beta}{2}}r^{\beta-1}\delta^2\big)
\end{equation}
$$
=C\big( A^2R^{\frac{1+3\beta}{1-\beta}}+R^{-1}\delta^2\big)=CR^{-1}\delta^2.
$$
\qed

\section{Regularity for Obstacle Problems (Proof of Main Theorem)}

\subsection{Heuristic arguments}
In this section we combine the results from our previous sections and
prove our main theorem, i.e. optimal regularity of the solution
to the obstacle problem with Hamilton-Jacobi obstacle.

The proof is again somewhat involved, so before we start let us
describe the idea.

Our goal is to prove optimal $C^{1,\alpha/2}$-regularity for minimizers
to the double obstacle problem with $h^\pm$ as obstacles. The proof
consists of several steps and an iteration argument. Here is a scheme of the steps:

\vspace{3mm}
\noindent
\textbf{Step 1)} From Corollary \ref{FirstRegObst} and Lemma \ref{PrelimObstReg} we already
have some regularity, in particular $C^{1,\frac{\alpha}{2+\alpha}}$-regularity for $u$ and one-sided $C^{1,\frac{\alpha}{2+\alpha}}$-estimates for $h^\pm$ at points where $h^\pm=u$.

\vspace{3mm}
\noindent
\textbf{Step 2)} Having $C^{1,\beta}$-estimates for $u$ and one-sided $C^{1,\beta}$-estimates for $h^\pm$, we can apply Proposition \ref{nondeg}
and Proposition \ref{orthogonal}, which will give us a slight, let us say
an $\epsilon$, gain
in the H\"older exponent for one-sided estimates in the $x'$-directions
for $h^\pm$. That is, $h^\pm$ satisfy one-sided $C^{1,\beta+\epsilon}$-estimates in the $x'$-directions.

\vspace{3mm}
\noindent
\textbf{Step 3)} Using that $u$ solves the obstacle problem together with one-sided $C^{1,\beta+\epsilon}$-estimates we can show that $u\in C^{1,\beta+\epsilon}$.
 In particular, we have gained an $\epsilon=\epsilon(\alpha,\beta)$
in regularity as compared to Corollary \ref{FirstRegObst} and Lemma \ref{PrelimObstReg}.

\vspace{3mm}
\noindent
\textbf{Step 4)} We can iterate Step 2) and 3) to gain more regularity of $u$, but in order
to fully utilize the $\epsilon$ gain in regularity we need
to be able to control the $C^{1,\beta+\epsilon}$-norm of $u$.
As it turns out there is a constant $\xi$ such that
if $[a]_{C^\alpha}\le \xi$ then we will get good control over
the $C^{1,\beta+\epsilon}$-norm of $u$. So we will rescale
$u$ and $a$ by a factor $\tau$, where $\tau$ is chosen such that
$[a(\tau x)]_{C^{\alpha}}\le \xi$. With this rescaling we can
iterate Step 2) and  3) to gain another $\epsilon$ in the H\"older exponent
\textsl{and also preserve the H\"older norm}.

Iterating Step 2) and  3) we will see that $u\in C^{1,\beta}$
for any $\beta< \alpha/2$, but \textsl{with a uniform bound on the
$C^{1,\beta}$-norm}. The Theorem follows.

In reality, the proof will be somewhat more involved as we have different
regularity on different scales which will result in some technical
issues.

With this strategy in mind let us turn to the proof of the Main Theorem.

\subsection{Proof of the Main Theorem}
 Without loss of generality we may assume that
$|\nabla u(0)|=|u(0)|=0$ ---$\nabla u$ exists by Corollary \ref{FirstRegObst}--- and that $a(x)=e_n+b(x)$ with $|b(x)|\le A|x|^\alpha$ in $B_1$.
If this is not true, we may subtract $u(0)+\nabla u(0)\cdot x$ from $u$ and $f$, add $\nabla u(0)$
to $a$ and rotate the coordinate system to obtain this situation.

We may also rescale $u_\tau(x)=\frac{u(\tau x)}{\tau}$ with $\tau=(\xi/A)^{1/\alpha}$
for a $\xi$ depending only on $\alpha$ and $n$ to be determined later. The scaled solution
$u_\tau$ will then minimize the Dirichlet energy in $B_{R_\tau}$ with $R_\tau=1/\tau$ and constraints
$h^\pm_\tau$ solving
$$
\pm|\nabla h^\pm_\tau-a_\tau|^2=\pm 1,
$$
where $a_\tau(x)=a(\tau x)$. We set $A_\tau:= [a_\tau]_{C^\alpha(B_{1/\tau})}= \xi$.

From Corollary \ref{FirstRegObst} we have $\osc_{B_r} u_\tau \le C A_\tau^{1/(2+\alpha)}r^{1+\alpha/(2+\alpha)}$ whenever $r\le A_\tau^{1/2}\min(R_\tau^{(2+\alpha)/2},R_\tau)$
and also ---using Lemma \ref{PrelimObstReg}---
$\sup_{B_r}h^+_\tau\le C A_\tau^{1/(2+\alpha)}r^{1+\alpha/(2+\alpha)}$.
That is enough to apply Lemma \ref{lemmaA} and Proposition \ref{orthogonal}
with $\beta_0=\alpha/(2+\alpha)$ and
$r\le \min(A_\tau^{1/2}R_\tau^{1/(1-\beta_0)},R_\tau)$.

Let us assume, for now, that $r\le R_\tau\le A_\tau^{1/2}R_\tau^{1/(1-\beta_0)}$. Then we have that for 
every $x_n\in(0,A_\tau^{-{1-\beta_0\over 2}}r^{1-\beta_0})$ and some
$p(x_n)=p$ in the super-differential of $h^+$ at $(0,x_n)$,
\begin{equation}\label{OneSidedforH1}
\sup_{x'\in B_r'}\big( h_\tau^+(x',x_n)-h_\tau^+(0,x_n)- p\cdot x' \big)\le
CA_\tau^{\frac{1}{2}-\frac{\alpha\beta_0}{2(2+\alpha-2\beta_0)}}r^{1+\frac{\alpha}{2+\alpha-2\beta_0}},
\end{equation}
if $r\le \textrm{min}\big( A_\tau^{\frac{6\beta_0+\alpha\beta_0-2-\alpha}{2(4\beta_0-\alpha)}},
A_\tau^{-\frac{\beta_0}{\alpha-2\beta_0}}, A_\tau^{\frac{\beta_0}{2}}R_\tau^{\frac{2+\alpha-2\beta_0}{2}}\big)$.

With 
$$
F(r)=CA_\tau^{\frac{1}{2}-\frac{\alpha\beta_0}{2(2+\alpha-2\beta_0)}}r^{1+\frac{\alpha}{2+\alpha-2\beta_0}},
$$
it follows from (\ref{OneSidedforH1}) that $u$ satisfies the assumptions of 
Proposition \ref{justneedwhatevername}, that is,
there is a constant 
$\tilde{C}$ such that
\begin{equation}\label{OneSidedforH12}
\textrm{osc}_{x\in B_r}u_\tau \le
\tilde{C}A_\tau^{\frac{1}{2}-\frac{\alpha\beta_0}{2(2+\alpha-2\beta_0)}}r^{1+\frac{\alpha}{2+\alpha-2\beta_0}},
\end{equation}
if $r\le \rho_\tau:= \textrm{min}\big( A_\tau^{1/2}R_\tau^{1/(1-\beta_0)},
A_\tau^{\frac{6\beta_0+\alpha\beta_0-2-\alpha}{2(4\beta_0-\alpha)}},
A_\tau^{-\frac{\beta_0}{\alpha-2\beta_0}},A_\tau^{\frac{\beta_0}{2}}R_\tau^{\frac{2+\alpha-2\beta_0}{2}}\big)$.

$ $

The estimate (\ref{OneSidedforH12}) gives us better regularity for small $r$.
We would want to iterate that statement. Let us rewrite the statement with $\beta_1=\alpha/(2+\alpha-2\beta_0)$:
\begin{equation}\label{herewenieedxi}
\textrm{osc}_{B_r}u_\tau \le \big(\tilde{C}A_\tau^{\frac{(1-\beta_0)\beta_1}{2}}\big) A_\tau^{\frac{1-\beta_1}{2}}r^{1+\beta_1}
\le A_\tau^{\frac{1-\beta_1}{2}}r^{1+\beta_1};
\end{equation}
here the last inequality is valid if $A_\tau$ is small enough, however the size of $A_\tau$ is independent
of $\beta_0$ and $\beta_1$ (given that $1/2\ge \beta_0,\beta_1\ge \alpha/(2+\alpha)$).
As $A_\tau=\xi$ by our rescaling,
let us choose $\xi$ as the largest
constant $\le 1$ such that the last inequality holds for all $1/2 \ge \beta_0,\; \beta_1\ge \alpha/(2+\alpha)$. 
Then $\xi$ depends only on $\alpha$ and $n$.

We are going to apply Proposition \ref{nondeg} and then Proposition \ref{orthogonal} again with $\beta=\beta_1$. It is  natural to split the proof
into the four cases
\begin{enumerate}
\item $R_\tau\ge 1 \textrm{ (in which case $A\ge \xi$)}$ and $R_\tau\le A_\tau R_\tau^{\frac{1}{1-\beta_0}}$,
\item $R_\tau\le 1 \textrm{ (in which case $A\le \xi$)}$ and $R_\tau\ge A_\tau R_\tau^{\frac{1}{1-\beta_0}}$,
\item $R_\tau\ge 1 \textrm{ (in which case $A\ge \xi$)}$ and $R_\tau\ge A_\tau R_\tau^{\frac{1}{1-\beta_0}}$,
\item $R_\tau\le 1 \textrm{ (in which case $A\le \xi$)}$ and $R_\tau\le A_\tau R_\tau^{\frac{1}{1-\beta_0}}$.
\end{enumerate}
Here (1) and (2) relate to the first and second block of values of $Q$
in Proposition \ref{orthogonal}, respectively. These cases are stable 
in the following sense: when $\beta$
increases in Case (1) or (2), then $A_\tau R_\tau^{\frac{1}{1-\beta}}$ increases if $R_\tau\ge 1$
and decreases if $R_\tau\le 1$.
That means that we stay in the same Case (1) or (2), respectively for larger
$\beta$.
In Case (3) and (4) this is no longer true, so
we might have to use one of the blocks of $Q$ for a finite number of iterations
and then switch over to the other block. 

We will finish the proof in Case 1 first.

\textbf{Case (1) ($R_\tau\ge 1$ and $R_\tau\le A^{1/2}_\tau R_\tau^{1/(1-\beta_0)}$):} In this case we may iterate and get a $\beta_2=\alpha/(2+\alpha-2\beta_1)$ such that (\ref{OneSidedforH12}) holds
with $\beta_2$ replacing $\beta_0$. We may iterate indefinitely
to obtain an increasing sequence of $\beta_j$ such that $\beta_j=\alpha/(2+\alpha-2\beta_{j-1})$,
and (\ref{OneSidedforH12}) holds with $\beta_j$ replacing $\beta_0$ for each $j\in \N$.
It is easy to see that $\beta_j\to \alpha/2$ as $j\to\infty$.
Using that $A_\tau=\xi$ is a constant depending only on $\alpha$ and $n$ and that $R_\tau\ge 1$, it follows that
\begin{equation}\label{forrgreaterthan1}
\textrm{osc}_{B_r}u_\tau\le \tilde{C}A_\tau^{\frac{4-\alpha^2}{8}}r^{1+\frac{\alpha}{2}}
\end{equation}
$$
\textrm{for }r\le \inf_{j}\Big(\textrm{min}\big( A_\tau^{\frac{6\beta_j+\alpha\beta_j-2-\alpha}{2(4\beta-\alpha)}},
A_\tau^{-\frac{\beta_j}{\alpha-2\beta_j}},A_\tau^{\frac{\beta_j}{2}}R_\tau^{\frac{2+\alpha-2\beta_j}{2}}\big)\Big)\le C_\alpha. 
$$ 

Rescaling back to $u$ we obtain from (\ref{forrgreaterthan1}) as well as 
the definition of $\tau$ that
\begin{equation}\label{rescaledrgreater1}
\textrm{osc}_{B_r} u = \tau \> \textrm{osc}_{B_{r/\tau}} u_\tau\le
C_1 \tau^{-\alpha/2}r^{1+\alpha/2}\le C_2 A_\tau^{1/2}r^{1+\alpha/2},
\end{equation}
for $r\le C_3 A_\tau^{-\frac{1}{\alpha}}$;
here $C_1,C_2$ and $C_3$ are constants 
depending only on $\alpha$ and $n$.

When $r\ge C_3 A_\tau^{-\frac{1}{\alpha}}$ we obtain from 
$\beta_0=\frac{\alpha}{2+\alpha}$ that
$$
\textrm{osc}_{B_r}u\le C_4 A_\tau^{\frac{1-\beta_0}{2}}r^{1+\beta_0}\le C_5 A_\tau^{\frac{1}{2}}r^{1+\frac{\alpha}{2}},
$$
where $C_4$ and $C_5$ are constants 
depending only on $\alpha$ and $n$. 
Combining the two estimates proves the Theorem in Case (1).

$ $

\textbf{Case (2) ($R_\tau\le 1$ and $R_\tau\ge A_{\tau}^{1/2}R_\tau^{1/(1-\beta_0)}$):}

When $A_\tau^{1/2}R_\tau^{\frac{1}{1-\beta_0}} =\inf(A_\tau^{1/2}R_\tau^{\frac{1}{1-\beta_0}},R_\tau)$, we
have to use the second block of values for $Q$ in Proposition \ref{orthogonal}.
We deduce that
$$
\sup_{B_r'(0,x_n)} |h^+(x)-h^+(0,x_n)-p\cdot x|\le C A^{\frac{1}{2}-\frac{\alpha\beta_0}{2(2+\alpha-2\beta_0)}}r^{1+\frac{\alpha}{2+\alpha-2\beta_0}}
$$
if $r\le\min\big( A^{\frac{6\beta_0+\alpha\beta_0-2-\alpha}{2(4\beta_0-\alpha)}},A^{-\frac{\beta_0}{\alpha-2\beta_0}}, A^{\frac{2+\alpha}{4}}R^{\frac{2+\alpha-2\beta_0}{2(1-\beta_0)}}\big)$.

Using Proposition \ref{justneedwhatevername} and iterating, we obtain that for
$\beta_j$ defined above,
\begin{equation}\label{Berkeley1}
\textrm{osc}_{x\in B_r}u_\tau\le
\tilde{C}A_\tau^{\frac{1}{2}-\frac{\alpha\beta_j}{2(2+\alpha-2\beta_j)}}r^{1+\frac{\alpha}{2+\alpha-2\beta_j}}
\end{equation}
for $r\le \textrm{min}\big(A_{\tau}^{\frac{6\beta_j+\alpha\beta_j-2-\alpha}{2(4\beta-\alpha)}},A_{\tau}^{-\frac{\beta_j}{\alpha-2\beta_j}},A_{\tau}^{\frac{2+\alpha}{4}}R_\tau^{\frac{2+\alpha-2\beta_j}{2(1-\beta_j)}}\big)$.

As in (\ref{herewenieedxi}) we conclude that
\begin{equation}\label{Berkeley127th}
\textrm{osc}_{x\in B_r}u_\tau\le
A_\tau^{\frac{1}{2}-\frac{\alpha\beta_j}{2(2+\alpha-2\beta_j)}}r^{1+\frac{\alpha}{2+\alpha-2\beta_j}}
\end{equation}
for $r\le \textrm{min}\big(A_{\tau}^{\frac{6\beta_j+\alpha\beta_j-2-\alpha}{2(4\beta-\alpha)}},A_{\tau}^{-\frac{\beta_j}{\alpha-2\beta_j}},A_{\tau}^{\frac{2+\alpha}{4}}R_\tau^{\frac{2+\alpha-2\beta_j}{2(1-\beta_j)}}\big)$.

Sending $j$ to infinity and using that $R_\tau<1$ 
and that $\frac{2+\alpha-2\beta_j}{2(1-\beta_j)}$
is decreasing in $\beta_j$ we obtain that
$$
\textrm{osc}_{x'\in B_r}u_\tau\le
C_6A_\tau^{\frac{1}{2}-\frac{\alpha^2}{8}}r^{1+\frac{\alpha}{2}}
$$
for $r\le A_\tau^{\frac{2+\alpha}{4}} R_\tau^{\frac{2}{2-\alpha}}$,
where $C_6$ is a constant 
depending only on $\alpha$ and $n$. 

Using Proposition \ref{justneedwhatevername} as well as 
the second case in the second block in Proposition \ref{orthogonal}
in iteration, we also see that
$$
\textrm{osc}_{x'\in B_r'}u_\tau\le C_7 R_\tau^{-1} r^2
$$
for $r\ge A_\tau^{\frac{2+\alpha}{4}} R_\tau^{\frac{2}{2-\alpha}}$,
where $C_7$ is a constant 
depending only on $\alpha$ and $n$.

Scaling back and using that $A_\tau$ is a constant depending only on $\alpha$ and $n$ we obtain that
$$
\textrm{osc} u \le 
\left\{
\begin{array}{ll}
C_8 A_\tau^{\frac{1}{2}}r^{1+\frac{\alpha}{2}} & \textrm{ if }r\le C_{10} A_\tau^{\frac{2}{2-\alpha}}, \\
C_9 r^2 & \textrm{ if }r\ge C_\alpha A_\tau^{\frac{2}{2-\alpha}},
\end{array}\right.
$$
where $C_8,C_9$ and $C_{10}$ are constants 
depending only on $\alpha$ and $n$. This concludes the proof in Case (2).

\textbf{Case (3) ($R_\tau\ge 1$ and $R_\tau\ge A_{\tau}^{1/2}R_\tau^{1/(1-\beta_0)}$):}
We start as in Case (2) and deduce that
for all $j$ such that 
\begin{equation}\label{star27th}
R_\tau\ge A_\tau^{\frac{1}{2}}R_\tau^{\frac{1}{1-\beta_{j-1}}},
\end{equation}
we have
$$
\sup_{B_r'(0,x_n)} |h^+(x)-h^+(0,x_n)-p\cdot x|\le C A^{\frac{1}{2}-\frac{\alpha\beta_j}{2(2+\alpha-2\beta_j)}}r^{1+\frac{\alpha}{2+\alpha-2\beta_j}}
$$
for  $r\le\min\big( A^{\frac{6\beta_j+\alpha\beta_j-2-\alpha}{2(4\beta_j-\alpha)}},A^{-\frac{\beta_j}{\alpha-2\beta_j}}, A^{\frac{2+\alpha}{4}}R^{\frac{2+\alpha-2\beta_j}{2(1-\beta_j)}}\big)$.

If (\ref{star27th}) holds for all $j$ then we are done as in Case (2).
If (\ref{star27th}) is not true then there is a largest $j_0$ such that 
the inequality holds for all $j\le j_0$. It follows that
$R_\tau\ge 1$ and $R_\tau\le A_\tau^{1/2}R_\tau^{\frac{1}{1-\beta_{j_0}}}$. 
These are the assumptions in the iteration in Case (1), so we may proceed as in 
Case (1) and obtain the statement of the Theorem.

\textbf{Case (4) ($R_\tau\le 1$ and $R_\tau\le A_{\tau}^{1/2}R_\tau^{1/(1-\beta_0)}$):}
We start as in Case (1) and deduce that as long as $R_\tau\le A_\tau^{1/2}R_\tau^{1/(1-\beta_{j-1})}$,
\begin{equation}
\textrm{osc}_{x\in B_r}u_\tau \le
A_\tau^{\frac{1}{2}-\frac{\alpha\beta_j}{2(2+\alpha-2\beta_j)}}r^{1+\frac{\alpha}{2+\alpha-2\beta_j}},
\end{equation}
for $r\le  \textrm{min}\big( A_\tau^{1/2}R_\tau^{1/(1-\beta_j)},
A_\tau^{\frac{6\beta_j+\alpha\beta_j-2-\alpha}{2(4\beta_j-\alpha)}},
A_\tau^{-\frac{\beta_j}{\alpha-2\beta_j}},A_\tau^{\frac{\beta_j}{2}}R_\tau^{\frac{2+\alpha-2\beta_j}{2}}\big)$.

If $R_\tau\le A_\tau^{1/2}R_\tau^{1/(1-\beta_{j-1})}$ for all $j$ then we are done.
If there is a $j_0$ such that $R_\tau\ge A_\tau^{1/2}R_\tau^{1/(1-\beta_{j_0})}$
then we are in the situation of Case (2), so we may proceed as in 
Case (2) and obtain the statement of the Theorem.

This finishes the proof. \qed

\section{Appendix: Remarks on Viscosity Solutions For Hamilton-Jacobi Equations.}\label{appendix}

In this appendix we will remind ourselves of some basic properties of viscosity
solutions for Hamilton-Jacobi equations. Most of the results in the
appendix can
be found in \cite{L}. However the exposition in \cite{L} is quite sketchy at
points and some of the results that we need are not explicitly proved. We will
therefore, for the readers convenience, supply some details, without claiming any 
originality. 
For the original and classical papers on the theory of viscosity solutions we 
direct the reader to \cite{CranEvanLion}, \cite{CranLions1} and 
\cite{CranLions2}.
First we will state the definition of viscosity
solutions for  first order Hamilton-Jacobi equations.
\begin{definition}
We say that $u\in C^0(\Omega)$ is a viscosity solution to
$$
H(x,u,\nabla u)=f(x) \textrm{ in } \Omega,
$$
if for every $x^0\in \Omega$ and $\phi\in C^1(B_r(x^0))$ such that
$u(x^0)=\phi(x^0)$ the following holds:
\begin{enumerate}
\item
if $u(x^0)=\phi(x^0)$ and $u(x)-\phi(x)$ has a local maximum at $x^0$, 
\linebreak
then $H(x^0,\phi(x^0),\nabla \phi(x^0))\le f(x^0)$,
\item if $u(x^0)=\phi(x^0)$ and
$u(x)-\phi(x)$ has a local minimum at $x^0$, \linebreak
then $H(x^0,\phi(x^0),\nabla \phi(x^0))\ge f(x^0)$.
\end{enumerate}
\end{definition}
This definition turns out to be the right one for first order Hamilton-Jacobi equations 
in the sense that it provides strong existence and uniqueness results.

Viscosity solutions can be obtained by the vanishing viscosity method.
That method also implies one-sided estimates for the second derivatives
of solutions to smooth convex first order Hamilton-Jacobi equations
which we will prove next.
\begin{lem}\label{OneSidedEstimates}
Let $u\in C^0(B_R\times(-R,0))$ be a viscosity solution to
$$
\frac{\partial u}{\partial t}+H(x,u,\nabla u)=m(x)
$$
such that $H\in C^2$ and $$p\cdot D^2_pH(x,a,p)\cdot p\ge c|p|^2 
\textrm{ for every } (x,a,p).$$
Moreover we assume that $m\in C^{1,1}$. Then the distributional second derivatives satisfy
$$
\frac{\partial^2 u}{\partial x_i^2}\le \frac{C_0}{R}
\textrm{ in } B_{R/2}\times (-R/2,0)
$$
for all spatial directions $x_i$. The constant $C_0$ depends only on $\|D^2 m \|_{L^\infty (B_R)}$, $\sup_{B_R}|\nabla u|$ and $H$.

A similar statement holds for the time independent case, that is if $u\in C^0(B_R)$ is a viscosity solution to
$$
H(x,u,\nabla u)=m(x)
$$
such that $H\in C^2$ and $$p\cdot D^2_pH(x,a,p)\cdot p\ge c|p|^2 
\textrm{ for every } (x,a,p)$$
and $m\in C^{1,1}$, then
$$
\frac{\partial^2 u}{\partial x_i^2}\le \frac{C_0}{R}
\textrm{ in } B_{R/2}
$$
for all spatial directions $x_i$. The constant $C_0$ depends only on $\|D^2 m \|_{L^\infty (B_R)}$, $\sup_{B_R}|\nabla u|$ and $H$.
\end{lem}
\textsl{Proof:} We will only prove the first statement for parabolic
Hamilton-Jacobi equations. The proof for elliptic Hamilton-Jacobi
equations is similar. Alternatively the elliptic proof may be derived by considering $t$ to be a dummy variable.

Moreover, as the proof is well known and included only for the sake of completeness, 
we will consider the slightly simpler case $H(x,a,p)=H(p)$. The general case
can be handled similarly.

Rescaling $u(Rx,Rt)/R$ if necessary we may assume that $R=1$.
Let $\chi\in C^\infty(B_1\times (-1,0))$ be a non negative function 
such that $\chi=1$ in
$B_{1/2}\times (-1/2,0)$ and $\chi=0$ close to $\{t=-1\}$ and
$\partial B_1\times (-1,1)$. Furthermore we may assume that
$|\nabla \chi|, |\Delta \chi|, |\nabla \chi|^2/|\chi|, |\partial \chi/\partial t|\le C$ in the support of $\chi$.

We are going to argue by the method of vanishing viscosity. That is,
we will approximate $\frac{\partial u}{\partial t}+H(\nabla u)=m(x)$
by the equation
\begin{equation}\label{vanishingviscosity}
-\epsilon \Delta u +\frac{\partial u}{\partial t}+H(\nabla u)=m(x)
\end{equation}
in order to deduce the desired estimate independent of $\epsilon>0$.
The Lemma follows from uniform convergence by letting $\epsilon\to 0$.

Differentiating equation (\ref{vanishingviscosity}) twice in the direction
$x_i$ we get
$$
-\epsilon \Delta u_{ii} +\frac{\partial u_{ii}}{\partial t}+H'_p(\nabla u)\cdot
\nabla u_{ii}+ \big(\nabla u_i \cdot H''_{pp}  \big)\cdot \nabla u_i=m_{ii}(x).
$$
Writing $w=\chi u_{ii}$ we see that
$$
-\epsilon \Delta w=-\chi\frac{\partial u_{ii}}{\partial t}-\chi
\big(\nabla u_i \cdot H''_{pp}  \big)\cdot \nabla u_i
$$
$$
-\chi
H'_p(\nabla u)\cdot \nabla u_{ii}+\chi m_{ii}-2\epsilon \nabla \chi\cdot
\nabla u_{ii}-\epsilon u_{ii}\Delta \chi.
$$

At a point $(x^0,t_0)$ where $w$ attains its positive  supremum we have
$-\Delta w\ge 0$, $\nabla w=0$ and $\partial w/\partial t\ge 0$.
The last two conditions are equivalent to
$$
\nabla u_{ii}(x^0,t_0)=-\frac{u_{ii}(x^0,t_0)\nabla \chi(x^0,t_0)}{\chi(x^0,t_0)},
$$
$$
\chi(x^0,t_0)\frac{\partial u_{ii}(x^0,t_0)}{\partial t}\ge
-u(x^0,t_0)\frac{\partial \chi(x^0,t_0)}{\partial t}.
$$
Using this together with the convexity assumption on $H$, we end up with
$$
0\le u_{ii}\frac{\partial \chi}{\partial t}-c\chi|\nabla u_i|^2+u_{ii} H'_p(\nabla u)\cdot\nabla \chi + \chi m_{ii}+2\epsilon \frac{|\nabla \chi|^2}{\chi}u_{ii}-
\epsilon u_{ii}\Delta \chi
$$
at $(x^0,t_0)$.
Rearranging terms we get
$$
c\chi|\nabla u_i|^2\le \big(\frac{\partial \chi}{\partial t}+H'_p(\nabla u)\cdot \nabla \chi+2\epsilon \frac{|\nabla \chi|^2}{\chi}-\epsilon\Delta \chi \big)u_{ii}+m_{ii}\chi.
$$
Using that ---by our choice of $\chi$ and assumptions on $H$--- the terms in the parenthesis
may be estimated by a constant $C$ depending only on $\sup_{B_R}|\nabla u|$ and $H$,
and observing that
the final term $m_{ii}\chi$ is also bounded by the assumption $m\in C^{1,1}$,
we obtain that at the point where $w$ attains its supremum,
$$
c\chi|\nabla u_i|^2\le C|u_{ii}|+C,
$$
where $C$ depends on
$\sup_{B_1\times (-1,0)}m_{ii}$, $\sup_{B_1\times (-1,0)}|\nabla u|$ and
$\sup |H'_p(\nabla u)|$. Multiplying both sides by $\chi$ implies that $|w|\le C$. Which in terms of $u_{ii}$ becomes
$$\sup_{B_{1/2}\times (-1/2,0)}u_{ii}\le \sup_{B_1\times (-1,0)}\chi u_{ii}\le C.$$\qed

Later on we will also need good stability estimates for solutions, proved in
the next two Lemmas.
\begin{lem}\label{ComparisonElliptic}
Let $H \in C^2$ and let $u\in C^{0,1}(\overline{B_R})$ and $v\in C^{0,1}(\overline{B_R})$ be
viscosity solutions in $B_R$ to
$$
H(x,\nabla u)=n(x)
$$
and
$$
H(x, \nabla v)=m(x)
$$
where $$p\cdot D^2_p H(x,p)\cdot p\ge c|p|^2\textrm{ for every } (x,p),$$
$m\ge n\ge \lambda >0$, $u\le v$ on $\partial B_R$ and
$\|u\|_{L^{\infty}(B_R)},\|v\|_{L^{\infty}(B_R)}\le K$. Then
$$
u\le v\le u+C(\lambda,K)R\max\big( \sup_{B_R}(m-n), \sup_{\partial B_R}(v-u)\big).
$$
\end{lem}
\textsl{Proof:} It is sufficient to prove the Lemma when $\sup_{\partial B_R}(m-n)$ and $\sup_{B_R}(v-u)$
are small. By a rescaling we may assume that $R=1$.

\noindent\textbf{Claim 1:} {\sl Let
$$
w_1=e^{K+1}-e^{-u}.
$$
Then there exists a function $F(x,a,p)$, convex in $p$, such that
\begin{equation}\label{uiop}
\frac{\partial F(x,a,p)}{\partial a}\ge \lambda,
\end{equation}
and in the viscosity sense
$$
F(x,w_1,\nabla w_1)=0.
$$
Moreover $w_2=e^{K+1}-e^{-v}$ is a viscosity solution to
$$
F(x,w_2,\nabla w_2)=\big(e^{K+1}-w_2 \big)(m-n).
$$}
\textsl{Proof of Claim 1:} A simple calculation shows that
$$
\big( e^{K+1}-w_1 \big)H(x, \nabla w_1/(e^{K+1}-w_1))=(e^{K+1}-w_1)m,
$$
so if we denote
$$
F(x,a,p)=(e^{K+1}-a)H\big( x,\frac{p}{e^{K+1}-a}\big)-(e^{K+1}-a)m,
$$
the first statement in the claim follows. That $F$ is convex in $p$ follows from
the convexity of $H$. A direct calculation yields
$$
\frac{\partial F(x,a,p)}{\partial a}=-H \big(x,\frac{p}{e^{K+1}-a}\big)+
H,_p \big(x,\frac{p}{e^{K+1}-a}\big)\cdot\frac{p}{e^{K+1}-a}+m\ge m\ge \lambda
$$
where we have used convexity of $H$. That
$$
F(x,w_2,\nabla w_2)=\big(e^{K+1}-w_2 \big)(m-n).
$$
follows from a simple calculation. The claim follows.

$ $

Next we observe that we may assume that $\Delta w_1\le C_\Delta$ for some
constant $C_\Delta$: we may regularize $F$ so it becomes $C^2$.
Denote the regularized version of $F$ by $F_\tau$ and let $w_1^\tau$
solve
$$
F_\tau(x,w_1^\tau,\nabla w_1^\tau)=0,
$$
then by Lemma \ref{OneSidedEstimates}, $\Delta w_1^\tau\le C_\tau$. We
also know that $w_1^\tau\in C^{\alpha}$ uniformly for all $\tau>0$, thus we
find a subsequence $\tau_j\to 0$ such that $w_1^{\tau_j}\to w_1^0$
uniformly. Using uniqueness results for viscosity solutions it is easy to see
that $w_1^0=w_1$. So if we can show the Lemma for each $w_1^{\tau_j}$
the result follows for $w_1$. Thus we may assume that $\Delta w_1\le C_\Delta$
as long as our final estimate does not depend on $C_\Delta$. We will make 
several more regularizations in what follows. In order
to simplify notation, we will by the just explained 
argument assume that $\Delta w_1\le C_\Delta$, and continue
to work with $w_1$ and not with $w_1^{\tau_j}$.

Next, we apply a standard convolution type regularization to $w_1$ and $w_2$, that is we denote
$$
w_i^\delta=\int w_i(y)\phi_\delta(x-y)dy,
$$
where $\phi_\delta$ is a standard mollifier. Then
$$
F(x,w_1^\delta,\nabla w_1^\delta)=h_1^\delta 
$$
$$
\textrm{ and }F(x,w_2^\delta,\nabla w_2^\delta)=\big( e^{K+1}-w_2 \big)(m-n)+h_2^\delta,
$$
where $h_i^\delta \to 0$ uniformly in $L^p$ for all $p\in [1,{+\infty})$ as
$\delta\to 0$. To do this regularization we need to assume
that $u$ and $v$ are solutions in $B_{1+\delta}$ or prove the result
in $B_{1-\delta}$. However as $\delta\to 0$ this will not make any difference so
we may as well ignore this slight complication.

Finally we will denote by $w_2^{\delta,\epsilon}$ the solution to
$$
-\epsilon \Delta w_2^{\delta,\epsilon}+F(x,w_2^{\delta,\epsilon},\nabla w_2^{\delta,\epsilon})=h_2^\delta,
w_2^{\delta,\epsilon}=w_2^\delta \textrm{ on }\partial B_1.$$
With these regularizations we end up with
$$
-\epsilon \Delta (w_2^{\delta,\epsilon}-w_1^{\delta})+
F(x,w_2^{\delta,\epsilon},\nabla w_2^{\delta,\epsilon})-
F(x,w_1^{\delta},\nabla w_1^{\delta})
$$
$$
\le \big(e^{K+1}-w_2\big)(m-n)+\epsilon C_\Delta+h_2^\delta-h_1^\delta.
$$

Next, in order to linearize we define
$$
B(x)=\int_0^1 F_{,p}(x,w_1^{\delta},t\nabla w_2^{\delta,\epsilon}-(1-t)\nabla w_1^{\delta})dt.
$$
Then, taking equation (\ref{uiop}) into consideration,
$$
-\epsilon \Delta (w_2^{\delta,\epsilon}-w_1^{\delta})+B(x)\cdot(\nabla w_2^{\delta,\epsilon}-\nabla w_1^{\delta})+\lambda(w_2^{\delta,\epsilon}-w_1^{\delta})
$$
$$\le 
\big(e^{K+1}-w_2\big)(m-n)+\epsilon C_\Delta+h_2^\delta-h_1^\delta.
$$
Now let $f^{\delta,\epsilon}$ be the solution to
$$
-\epsilon\Delta f^{\delta,\epsilon}+B(x)\cdot \nabla f^{\delta,\epsilon}+\lambda f^{\delta,\epsilon}=\big(e^{K+1}-w_2 \big)(m-n) +C_\Delta\epsilon
$$
with boundary values $f^{\delta,\epsilon}=\sup_{\partial B_1}(w_2^{\delta,\epsilon}-w_1^{\delta,\epsilon})$
and $g^{\delta,\epsilon}$ be the solution to
$$
-\epsilon\Delta g^{\delta,\epsilon}+B(x)\cdot \nabla g^{\delta,\epsilon}+\lambda g^{\delta,\epsilon}=h_2^\delta-h_1^\delta,
$$
with boundary values $g^{\delta,\epsilon}=0$ on $\partial B_1$.
By the comparison principle
\begin{equation}\label{w2-w1lef+g}
w_2^{\delta,\epsilon}-w_1^{\delta,\epsilon}\le f^{\delta,\epsilon}+g^{\delta,\epsilon}.
\end{equation}

Since $h_2^\delta-h_1^\delta\to 0$ in $L^p$ for any $p<{+\infty}$ as $\delta\to 0$ we see that
$\sup_{B_1}\vert g^{\delta,\epsilon}\vert \to 0$ as $\delta\to 0$. Next we notice that where
$f^{\delta,\epsilon}$ attains its supremum we have $|\nabla f^{\delta,\epsilon}|=0$
and $\Delta f^{\delta,\epsilon}\le 0$. It follows that
$$
\sup_{B_1}f^{\delta,\epsilon}\le \max\Big( \sup_{\partial B_1}f^{\delta,\epsilon}, \sup_{B_1}\frac{1}{\lambda}\big(e^{K+1}-w_2 \big)(m-n) +\frac{C\epsilon}{\lambda} \Big).
$$
Letting first $\delta\to 0$ and then $\epsilon\to 0$ we see that equation
(\ref{w2-w1lef+g}) implies that
$$
w_2-w_1\le
\max\Big( \sup_{\partial B_1}(w_2-w_1), \sup_{B_1}\frac{1}{\lambda}\big(e^{K+1}-w_2 \big)(m-n)\Big).
$$
Writing this last inequality in terms of $u$ and $v$ we get
$$
e^{-u}-e^{-v}\le
\max\Big( \sup_{\partial B_1}\big(e^{-u}-e^{-v}\big), \sup_{B_1}\frac{e^{-v}}{\lambda}(m-n)\Big).
$$
Thus
$$
1-e^{u-v}\le e^{2K}\max\Big( \sup_{\partial B_1}\big(1-e^{u-v}\big), \sup_{B_1}\frac{1}{\lambda}(m-n)\Big).
$$
As pointed out in the beginning of the proof it is enough to show the
Lemma when $\sup_{B_1}(m-n)$ and $\sup_{\partial B_1}(v-u)$ are small.
Using
$\xi/2\le 1-e^{-\xi}\le 2\xi$ for small $\xi$, we end up with
$$
\sup_{B_1}(v-u)\le 2e^{2K}\max\big(\frac{1}{\lambda}\sup_{B_1}(m-n),2\sup_{\partial B_1}(v-u) \big).
$$
\qed

Next we need a comparison estimate for parabolic Hamilton-Jacobi
equations. The proof is very similar to the proof of Lemma
\ref{ComparisonElliptic} so we will omit it here with a reference to
\cite{L}.

\begin{lem}\label{ComparisonParabolic}
Let $H \in C^2$ and let
$u\in C^{0,1}(\overline{B_R}\times (-R,0))$ and
$v\in C^{0,1}(\overline{B_R}\times (-R,0))$ be
viscosity solutions in $B_R\times (-R,0)$ to
$$
\frac{\partial u}{\partial t}+H(\nabla u)=m(x)
$$
and
$$
\frac{\partial v}{\partial t}+H(\nabla v)=n(x),
$$
where $$p\cdot D^2H(p)\cdot p\ge c|p|^2 \textrm{ for every } p,$$
$n>m$, $u=v$ on
$\partial B_R$ for $t\in (-R,0)$ and $u(x,-R)=v(x,-R)$ for $x\in B_R$. Then
$$
u\le v\le u+CR\sup_{B_R\times (-R^2,0)}(n-m).
$$
\end{lem}
\textsl{Proof:} For a discussion of a proof see \cite{L}.\qed

From the regularity results in Lemma \ref{OneSidedEstimates} and 
\ref{ComparisonElliptic}  we can easily deduce
some elementary, non-optimal, one-sided estimates for solutions to
Hamilton-Jacobi equations, even when the data is not $C^2$. This incidentally 
provides our starting regularity for a bootstrap argument which in turn yields optimal
regularity.
\begin{lem}\label{PrelimObstReg}
Let $h$ be a viscosity solution to $|\nabla h-a|^2=1$ in $B_R$,
$a\in C^{\alpha}(\overline{B_R};\mathbb{R}^n)$ and
$[a(x)]_{C^\alpha(\overline{B_R})}\le A.$
Then for any $x^0\in B_{R/2}$ the super-differential of $h$ at $x^0$ is not 
empty, and for any $p$ in the super-differential of $h$ at $x^0$ we have
\begin{itemize}
\item[(1)]\label{toolongago}
$$
\sup_{B_r(x^0)}\big( h(x)-p\cdot(x-x^0)-h(x^0)\big) \le C A^{1/(2+\alpha)}r^{1+\alpha/(2+\alpha)}
$$
for $r\le \min(A^{1/2}R^{(2+\alpha)/2},R)$,
\item[(2)]
and
$$
\sup_{B_r(x^0)}\big( h(x)-p\cdot(x-x^0)-h(x^0)\big)\le C \frac{r^2}{R}
$$
for $A^{-1/\alpha}\ge r\ge A^{1/2}R^{(2+\alpha)/2}$.
\end{itemize}
The constant $C$ depends only on $n$.
\end{lem}
\textsl{Proof:} We may assume that $h(0)=0$.

First we notice that ---as $r\le R$--- if $r\ge A^{-1/\alpha}$ then  
$r\le A^{1/2}R^{(2+\alpha)/2}$.
Also, if
$r\ge A^{-1/\alpha}$ then 
$$
1= \sup_{B_r}|\nabla h-a|\ge \sup_{B_r}|\nabla h|-
\sup_{B_r}|a|=\sup_{B_r}|\nabla h|-Ar^\alpha.
$$
Using that $Ar^{\alpha}\ge 1$, it follows that
$$
\sup_{B_r}|h|\le r \sup_{B_r}|\nabla h| \le r\big(1 +Ar^\alpha \big)\le 
C A^{1/(2+\alpha)}r^{1+\alpha/(2+\alpha)}.
$$
This shows that if $r\ge A^{-1/\alpha}$ then both the assumption and the
conclusion in (1) hold.

We therefore only have to show the Lemma when $r\le A^{-1/\alpha}$.
We do this by a barrier type argument. We may assume that
$a(0)=0$. We also define the barrier $w$ as the solution to
$$
\begin{array}{ll}
|\nabla w|^2=1+(2 A \delta^{\alpha}+A^2\delta^{2\alpha}) & \textrm{in } B_{\delta}(x^0), \\
w=h & \textrm{on } \partial B_{\delta}(x^0),
\end{array}
$$
where $\delta$ is to be determined later.
Then $|\nabla w|\le \sqrt{1+(2 A \delta^{\alpha}+A^2\delta^{2\alpha})}\le 1+(2A\delta^\alpha+A^2\delta^{2\alpha})$ and thus
$$
|\nabla w-a|\ge|\nabla w|-|a|\ge\sqrt{1+(2 A \delta^{\alpha}+A^2\delta^{2\alpha})}-A\delta^\alpha 
$$
$$
= 1+\big(\sqrt{1+(2 A \delta^{\alpha}+A^2\delta^{2\alpha})}-(1+A\delta^\alpha) \big)=1, 
$$
so $w$ is a super-solution to $|\nabla h-a|^2=1$. Similarly,
$$
|\nabla w-a|^2\le 1+CA\delta^\alpha \big( 1+A\delta^\alpha \big).
$$
We may thus use the comparison
estimate for Hamilton-Jacobi equations (Lemma \ref{ComparisonElliptic}) 
and derive
$$
h\le w\le w(x^0)+p\cdot (x-x^0) +\frac{C|x-x^0|^2}{\delta}
$$
$$\le
h(x^0)+p\cdot (x-x^0) +\frac{C|x-x^0|^2}{\delta}+CA\delta^{1+\alpha}\big(1+A\delta^\alpha\big)
$$
for $p$ in the super-differential of $w$. Rearranging the terms
and taking the supremum in $B_r(x^0)$ for some $r\le\delta$ yields
\begin{equation}\label{anothercaseforbalance}
\sup_{B_r(x^0)}\big( h(x)-p\cdot(x-x^0)-h(x^0)\big) \le C \Big( \frac{r^2}{\delta}+A\delta^{1+\alpha}\big(1+A\delta^\alpha\big) \Big).
\end{equation}
We want to find the right balance between $r$ and $\delta$ that optimises (\ref{anothercaseforbalance}).
It is convenient to divide the last part of the proof into two cases:

{\bf Case 1:}  When $r\le \min(A^{1/2}R^{\frac{2+\alpha}{2}}, A^{-\frac{1}{\alpha}})$,
we use $\delta=A^{-\frac{1}{2+\alpha}}r^{\frac{2}{2+\alpha}}$ in (\ref{anothercaseforbalance})
and deduce that ---using $r\le A^{-1/\alpha}$ which implies that $A\delta^\alpha= A^{\frac{2}{2+\alpha}}r^{\frac{2\alpha}{2+\alpha}}\le 1$---
$$
\sup_{B_r(x^0)}\big( h(x)-p\cdot(x-x^0)-h(x^0)\big) 
$$
$$\le
C \Big( A^{\frac{1}{2+\alpha}}r^{1+\frac{\alpha}{2+\alpha}} +A\delta^{1+\alpha}\big(1+A\delta^\alpha\big) \Big)\le
CA^{\frac{1}{2+\alpha}}r^{1+\frac{\alpha}{2+\alpha}}.
$$
From the definition of $r$ and $\delta$ it is easy to check that $r\le \delta\le R$.

{\bf Case 2:} Next, if $A^{1/2}R^{\frac{2+\alpha}{2}}\le r\le A^{-1/\alpha}$ then we use $\delta=R$ in (\ref{anothercaseforbalance}) and deduce that ---using
$A^{1/2}R^{\frac{2+\alpha}{2}}\le r$,
$A^{1/2}R^{\frac{2+\alpha}{2}}\le r\le A^{-1/\alpha}$ and $AR^\alpha\le 1$---
$$
\sup_{B_r(x^0)}\big( h(x)-p\cdot(x-x^0)-h(x^0)\big) \le C \Big( \frac{r^2}{R}+AR^{1+\alpha}\big(1+AR^\alpha\big) \Big)\le C\frac{r^2}{R}.
$$
Observe that this interval is empty unless $R\le A^{-1/\alpha}$. %That is if
%$R\le A^{-1/\alpha}$ then the Lemma follows from 1) and 2). 
\qed

We end this appendix with a lemma reminiscent of a Whitney extension theorem, which will be useful in the main text.

\begin{lem}\label{lemmaA}
Let $h\in C^{0,1}(\overline{B_1})$ and assume that $h$ satisfies the one-sided $C^{1,\beta}$-estimate
\begin{equation}\label{Hexpl}
h(x^0+x)\le h(x^0)+ p_{x^0}\cdot (x-x^0) +C_0|x|^{1+\beta}
\end{equation}
for every $p_{x^0}$ in the super-differential of $h$ at $x^0$ and every $x^0\in B_1$.
Moreover assume that $h(x)\ge -C_1|x|^{1+\beta}$ and that $h(0)=0$. Then
\begin{equation}\label{Hexpl2}
|p_{x}|\le E(\beta)(C_0+C_1)|x|^{\beta} \textrm{ for } x\in B_{1/2};
\end{equation}
here $E(\beta)$ depends continuously on $\beta$ for $\beta \in (0,1)$. 
\end{lem}
\textsl{Proof:} Take $y\in B_r$ and let $p_y=p$ be in the super-differential
of $h$ at $y$. As $h$ is semi-concave by (\ref{Hexpl}), we know that the 
super-differential of $h$ is non-empty for every $x^0$.
Notice that, due to (\ref{Hexpl2}), the super-differential of $h$ at the origin 
contains only
the zero vector.

Therefore $h(y)\le C_0r^{1+\beta}$.
Now consider $z=y-\epsilon p$ with
$$
\epsilon=\big( \delta |p|^{1-\beta} \big)^{\frac{1}{\beta}}
$$
for some small constant $\delta$. By one-sided estimates from above and below we have
$$
-C_\beta C_1\big( r^{1+\beta}+\epsilon^{1+\beta}|p|^{1+\beta} \big)
$$
$$\le 
-C_1|y-\epsilon p|^{1+\beta}\le h(z)\le h(y)-\epsilon |p|^2+
C_0\epsilon^{1+\beta}|p|^{1+\beta},
$$
where $C_\beta$ is a constant depending only on $\beta$ and $n$. Reordering terms and using that $h(y)\le C_0 r^{1+\beta}$ 
yields
$$
\big( \delta^{1/\beta}-(C_\beta C_1 +C_0)\delta^{(1+\beta)/\beta} \big)|p|^{\frac{1+\beta}{\beta}}\le (C_0+C_\beta C_1)r^{1+\beta}.
$$
Choosing $\delta=(2C_0+2C_\beta C_1)^{-1}$, the previous inequality becomes, for some $C$ depending only on $\beta$,
$$
|p|^{(1+\beta)/\beta}\le C(C_1+C_0)^{(1+\beta)/\beta}r^{1+\beta}.
$$
Taking both sides to the power of $\beta/(1+\beta)$, the lemma follows.\qed

\bibliographystyle{plain}
\bibliography{TorsionFinal.bib}

\end{document}